\documentclass[12pt,a4paper]{article}
\usepackage[dvips]{graphicx}
\usepackage{latexsym}
\usepackage{amssymb}
\usepackage{amsmath}
\usepackage{amscd}
\usepackage{mathrsfs}
\newtheorem{Th}{Theorem}[section]
\newtheorem{Co}[Th]{Corollary}
\newtheorem{Lem}[Th]{Lemma}

\newtheorem{Pro}[Th]{Proposition}
\newcommand{\demo}{\par\noindent{\it Proof. \/}\ }
\newcommand{\enD}{\hfill $\Box$\vspace{3truemm} \par}
\newcommand{\bx}{\mbox{\boldmath $x$}}
\newcommand{\bX}{\mbox{\boldmath $X$}}

\newcommand{\be}{\mbox{\boldmath $e$}}
\newcommand{\ba}{\mbox{\boldmath $a$}}
\newcommand{\bb}{\mbox{\boldmath $b$}}
\newcommand{\bv}{\mbox{\boldmath $v$}}
\newcommand{\by}{\mbox{\boldmath $y$}}
\newcommand{\bo}{\mbox{\boldmath $0$}}

\newcommand{\bxi}{\mbox{\boldmath $\xi$}}

\newcommand{\blambda}{\mbox{\boldmath $\lambda$}}
\newcommand{\bsigma}{\mbox{\boldmath $\sigma$}}
\newcommand{\bgamma}{\mbox{\boldmath $\gamma$}}
\newcommand{\bt}{\mbox{\boldmath $t$}}
\newcommand{\bn}{\mbox{\boldmath $n$}}

\newcommand{\sbxi}{\mbox{\scriptsize \boldmath$\xi$}}
\newcommand{\sbn}{\mbox{\scriptsize \boldmath$n$}}
\newcommand{\sbgamma}{\mbox{\scriptsize \boldmath$\gamma$}}

\newcommand{\sblambda}{\mbox{\scriptsize \boldmath$\lambda$}}

\newcommand{\R}{{\mathbb R}}
\newcommand{\lon}{\longrightarrow}
\begin{document}
\title{Lightlike hypersurfaces along spacelike submanifolds in anti-de Sitter space
 }
\author{Shyuichi IZUMIYA}

\date{\today}
\maketitle

%\bigskip
\begin{abstract}
Anti-de Sitter space is the Lorentzian space form with negative curvature.
In this paper we consider lightlike hypersurfaces along spacelike submanifolds in anti-de Sitter space with
general codimension.
In particular, we investigate the singularities of lightlike hypersurfaces as an application of
the theory of Legendrian singularities.
\end{abstract}
\renewcommand{\thefootnote}{\fnsymbol{footnote}}

\section{Introduction}
\label{intro}

Anti-de Sitter space is one of the Lorentz space forms which has rich geometric properties.
It is defined as a pseudo-sphere with negative curvature in semi-Euclidean space with index 2 which admits the biggest symmetry in  Riemannian or Lorentz space forms.
Anti-de Sitter space can be naturally considered as a Lorentzian version (generalization) of Hyperbolic space. Recently we discovered interesting geometric properties of submanifolds in Hyperbolic space as an
application of the theory of Legendrian/Lagrangian singularities\cite{Izu2,Izu11,IPT04,IPRT05}.
Therefore anti-de Sitter space might have rich geometric properties comparing with Hyperbolic space.
This is one of the motivations for the investigation of submanifolds in anti-de Sitter space from a mathematical viewpoint.
\par
On the other hand, anti-de Sitter space plays important roles in theoretical physics such as the theory of general relativity, the 
string theory and the brane world scenario etc \cite{BR02,M98,RS99,W98}.
It is well known that Lorentzian space forms are classified into three types depending on the value of the scalar curvature.
One of them is Lorentz-Minkowski space which has zero curvature.
The Lorentz space form with positive curvature is de Sitter space.
Anti-de Sitter space is a Lorentzian space form with negative curvature.
Recently, submanifolds in Lorentz-Minkowski space or de Sitter space have been well investigated (cf., \cite{Izu3,Izu4,IzuSM,IzuAG,Kosso1,Kosso2}).
However, there are not so many results on submanifolds in anti-de Sitter space, in particular  from the
viewpoint of singularity theory.
The lightlike hypersurfaces (i.e. the light-sheets in physics) are important objects in theoretical physics because they provide good models for  different types of horizons \cite{Ch, MTW}.
A lightlike hypersurface is generally a ruled hypersurface along a spacelike submanifold with codimension two whose rulings are lightlike geodesics. In this paper we consider lightlike hypersurfaces along spacelike submanifolds with general codimension in anti-de Sitter space.
Moreover, the lightlike hypersurface in anti-de Sitter $5$-space is a very important subject in the brane world scenario\cite{KR01,BR02,RS99}.
\par
In the meantime, tools in the theory of singularities
have proven to be useful 
description of geometrical properties of submanifolds immersed in
different ambient spaces, from both the local and global viewpoint
\cite{Izu2,Izu3,IPRT05,IzuSM,IzuAG,Little}.
The natural connection between geometry and singularities relies on
the basic fact that the contacts of a submanifold with the models of
the ambient space can be described by means of the analysis of the
singularities of appropriate families of contact functions, or
equivalently, of their associated  Legendrian maps
(\cite{Arnold1,Montaldi,Zak}).
When working in a Lorentz space form, the properties associated to the
contacts of a given submanifold with  lightcones have a
special relevance. 
In \cite{Izu3,IzuSM, Kasedou2}, it was constructed a framework for the study of
spacelike submanifolds with codimension two in Lorentz-Minkowski space or de Sitter space and discovered a Lorentz invariant 
concerning their contacts with models related to lightlike hyperplanes. 
The geometry described in this framework is called the {\it lightlike geometry} of spacelike submanifolds
with codimension two.
By using the invariants of lightlike geometry,
the singularities of lightlike hypersurfaces along spacelike submanifolds with codimension two in Lorentz-Minkowski space 
or de Sitter space
were investigated in \cite{Izu4,Izu08,Kasedou2}.
However, the situation is rather complicated for the general codimensional case.
The main difference from the Euclidean space (or, Hyperbolic space) case is the fiber of the canal hypersurface of a spacelike submanifold is neither connected nor
compact. 
In order to avoid the above difficulty, we arbitrarily choose a timelike future directed unit normal vector field
along the spacelike submanifold which always exists for an orientable manifold (cf., \S 3).
Then we construct the unit spherical normal bundle relative to the above timelike unit normal vector field, which can be considered as a codimension two spacelike canal submanifold of
the ambient space.
Therefore, we can apply the idea of the lightlike geometry of spacelike submanifolds of
the ambient space with codimension two. 
In this paper we apply the idea of this framework and the theory of Legendrian singularities to
investigate the singularities of lightlike hypersurfaces along spacelike submanifolds in anti-de Sitter space with general codimension.
\par
In \S 3 we construct the framework of the lightlike geometry of spacelike submanifolds with
general codimension analogous to \cite{IPRT05}.
The notion of lightlike hypersurfaces along spacelike submanifolds is introduced and  the basic properties are investigated in \S 4.
 The notion of the anti-de Sitter height functions families
is useful for the study of lightlike hypersurfaces (cf., \S 4).
The critical value set of the lightlike hypersurface along a spacelike submanifold is called
the {\it lightlike focal set} of the submanifold.
In \S 5 we show that the lightlike focal set of a spacelike submanifold is a point
if and only if the lightlike hypersurface along the submanifold is a subset of a lightcone (Proposition 5.1). Therefore, an anti-de Sitter lightcone is 
a model hypersurface of lightlike hypersurfaces.
The geometric meaning of the singularities of lightlike hypersurface is described by
the theory of contact of submanifolds with model hypersurfaces.
Moreover, as an application of the theory of Legendrian singularities, we show that
two lightlike hypersurfaces are locally diffeomorphic if and only if the types of the
contact of spacelike submanifolds with lightcones are the same in the sense of Montaldi\cite{Montaldi} under some generic conditions (Theorem 5.5).
In \S 6 we describe the case of codimension two as a special case.
We describe the detailed properties of lightlike focal sets of spacelike surfaces in 
anti-de Sitter $4$-space.
We also investigate spacelike curves in anti-de Sitter $4$-space as
the simplest case of a higher codimension in \S 7.
 \par
We shall assume throughout the whole paper that all the
maps and manifolds are $C^{\infty}$ unless the contrary is explicitly stated.

\section{Basic facts and notations on semi-Euclidean space with index 2}
\label{sec:1}
\par
In this section we prepare the basic notions on semi-Euclidean
(n+2)-space with index 2. For details of semi-Euclidean geometry,
see \cite{Oneil}.
Let $\R^{n+2}=\{(x_{-1}, x_0, x_1,\cdots,x_{n})|x_i \in \R\
(i=-1, 0,\cdots,n)\ \}$ be an (n+2)-dimensional vector space. For
any vectors $\bx=(x_{-1}, x_0, x_1,\cdots,x_{n})$ and $\by
=(y_{-1}, y_0, y_1,\cdots,y_{n})$ in ${\mathbb R}^{n+2},$ the {\it
pseudo scalar product \/} of $\bx$ and $\by$ is defined to be
$\langle\bx,\by\rangle =-x_{-1} y_{-1}-x_0 y_0 +\sum_{i=1}^{n}x_i
y_i$. We call $(\Bbb R^{n+2}, \langle ,\rangle )$ {\it
semi-Euclidean\/} (n+2)-{\it space with index 2\/} and write $\Bbb
R^{n+2}_2$ instead of $(\Bbb R^{n+2},\langle ,\rangle )$.
We say that a non-zero vector $\bx$ in $\Bbb R^{n+2}_2$ is {\it
spacelike\/}, {\it null\/} or {\it timelike\/} if
$\langle\bx,\bx\rangle>0,\langle\bx,\bx\rangle =0$ or
$\langle\bx,\bx\rangle <0$, respectively. The norm of the vector $\bx
\in \Bbb R^{n+2}_2$ is defined to be $\|\bx\|=\sqrt{|\langle\bx,
\bx\rangle|}$.
We define the signature of $\bx$ by
\[{\rm sign} (\bx)=\left\{
  \begin{array}{ccc}
  1\qquad\quad $\bx$\ \mbox{is\ spacelike,}\\

  \mbox{}0\qquad\quad \bx\ \mbox{is\ null,} \hspace*{\fill}\\

  -1\qquad \bx\ \mbox{is\ timelike.}
  \end{array}\right.
  \]
 For a non-zero vector $\bn\in \Bbb R^{n+2}_2$ and a real number $c$, we define the
 {\it hyperplane with pseudo-normal \/}$\bn$ by
$$
HP(\bn,c)=\{\bx \in \Bbb R^{n+2}_2 |\langle\bx,\bn\rangle=c\}.
$$
We call $HP(\bn,c)$ a {\it Lorentz hyperplane\/}, a {\it
semi-Euclidean hyperplane with index 2\/} or a {\it null
hyperplane\/} if $\bn$ is {\it timelike, spacelike or null
\/}respectively.

We now define an {\it anti-de Sitter $(n+1)$-space \/} (briefly, {\it AdS
$(n+1)$-space\/}) by
$$
AdS^{n+1}=\{\bx\in \Bbb R_2^{n+2}\ |\ \langle\bx,\bx\rangle =-1\}=H^{n+1}_1,
$$
a {\it unit pseudo $(n+1)$-sphere with index 2\/} by
$$
S^{n+1}_2=\{\bx\in \R_2^{n+2}\ |\ \langle\bx,\bx\rangle =1\}
$$
and a {\it {\rm (}closed{\rm )} nullcone\/} with vertex $\ba$ by\\
$$\Lambda_{a}^{n+1}=\{\bx\in \R_{2}^{n+2}|\langle\bx-\ba, \bx-\ba\rangle=0\}.
$$
In particular we denote that $\Lambda ^*=\Lambda ^{n+1}_0\setminus \{\bo\}$ and also call it
the {\it {\rm (}open{\rm )} nullcone\/}. Our main subject in this paper is $AdS^{n+1}$.
Since there are timelike closed curves in $AdS^{n+1}$, the causality of $AdS^{n+1}$ is violated.
In order to avoid such a situation, it is usually considered the universal covering
space $\widetilde{AdS}^{n+1}$ of $AdS^{n+1}$ in physics which is called the {\it universal Anti de Sitter space\/}.
We remark that the local structure of these spaces are the same.
\par
For any $\bx_1,\cdots, \bx_{n} \in \R^{n+2}_2$. We define a
vector $\bx_1\wedge \cdots\wedge \bx_n$ by
$$
\bx_1\wedge\cdots\wedge \bx_n= \vmatrix
-\be_{-1}&-\be_0&\be_1&\cdots&\be_{n}\vspace{2mm}\\
x^1_{-1}&x^1_0&x^1_1&\cdots&x^1_{n}\\
\vdots&\vdots&\vdots&\vdots&\vdots\\
x^{n}_{-1}&x^{n}_{0}&x^{n}_1&\cdots&x^{n}_{n}
\endvmatrix,
$$
where $\{\be_{-1}, \be_0, \be_1,\cdots,\be_{n}\}$ is the canonical
basis of $\Bbb R^{n+2}_2$ and $\bx_i=(x^i_{-1}, x^i_0,
x^i_1,\cdots,x^i_{n})$. We can easily check that
$$
\langle\bx,\ \bx_1\wedge\cdots\wedge
\bx_{n}\rangle=\textrm{det}(\bx, \bx_1,\cdots,\bx_{n}),
$$
so that $\bx_1\wedge\cdots\wedge \bx_{n}$ is pseudo-orthogonal to
any $\bx_i$\ (for\ $i=1,\cdots,n$).

\section{Spacelike submanifolds in anti-de Sitter space}
\label{sec:2}
 We introduce in this section the basic geometrical
tools for the study of spacelike submanifolds in
the anti-de Sitter $(n+1)$-space. 
\par
 Consider the orientation of $\R^{n+2}_2$
provided by the condition that $\textrm{det}(\be_{-1}, \be_0,\be _1,\cdots,\be_{n})>0.$
This orientation induces the orientation of $x_{-1}x_0$-plane,
so that it gives a time
orientation on $AdS^{n+1}$. 
If we consider the universal Anti de Sitter space $\widetilde{AdS}^{n+1},$ we can determine the
future direction.
\par
We consider a spacelike embedding $\bX:U\rightarrow AdS^{n+1}$ from an open subset $U\subset \R^s$
with $s+k=n+1.$ We write
$M=\bX(U)$ and identify $M$ and $U$ through the embedding $\bX.$ We
say that $\bX$ is {\it spacelike} if the tangent space $T_p  M$
consists only spacelike vectors (i.e. spacelike subspace) for any point $p\in M$. In this case, the pseudo-normal
space $N_p(M)$ in $\R^{n+2}_2$ is a $k+1$-dimensional semi-Euclidean space with index $2$
(cf. \cite{Oneil}). We write $N(M)$ as the pseudo-normal bundle in $\R^{n+2}_2$
over $M.$ 
On the pseudo-normal space $N_p(M),$ we have two kinds of pseudo spheres:
\begin{eqnarray*}
N_p(M;-1)& = & \{\bv\in N_p(M)\ |\ \langle \bv,\bv\rangle =-1\ \} \\
N_p(M;1)&= & \{\bv\in N_p(M)\ |\ \langle \bv,\bv\rangle =1\ \},
\end{eqnarray*}
so that we have two unit spherical normal bundles over $M$:
\[
N(M;-1)=\bigcup _{p\in M} N_p(M;-1)\ \mbox{and}\  N(M;1)=\bigcup _{p\in M} N_p(M;1).
\]
Then we have the Whitney sum decomposition
\[
T\R^{n+2}_2|M=TM\oplus N(M).
\]
By definition $\bX(u)$ is one of the timelike unit normal vectors of $M$ at $p=\bX(u),$ so that
$\bX \in N_p(M).$ Since $AdS^{n+1}$ is time oriented, we can arbitrarily choose an
adopted unit timelike normal section $\bn ^T(u)\in N_p(M)$ pseudo-orthogonal to
$\bX(u)$ even globally. Here, we say that $\bn^T$ is {\it adopted}
if 
\[
\textrm{det}(\bX(u),\bn^T(u),\be _1,\dots ,\be _{n})>0.
\]
Therefore we have the pseudo-orthonormal complement
$(\langle \bX(u),\bn ^T(u)\rangle _\R)^\perp$ in $N_p(M)$
which is a $(k-1)$-dimensional subspace of $N_p(M).$
We define a $(k-2)$-dimensional spacelike unit sphere in $N_p(M)$ by
\[
N_1(M)_p[\bn ^T]=\{\bxi \in N_p(M;1)\ |\ \langle \bxi, \bn ^T(p)\rangle =\langle \bxi,\bX(u)\rangle=0,p=\bX(u)\ \}.
\]
Then we have a {\it spacelike unit $k-2$-spherical bundle over $M$ with respect to $\bn ^T$} defined by
\[
N_1(M)[\bn ^T]=\bigcup _{p\in M} N_1(M)_p[\bn ^T].
\]
Since we have
$T_{(p,\xi)}N_1(M)[\bn^T]=T_pM\times T_\xi N_1(M)_p[\bn ^T],$
we have the canonical Riemannian metric on $N_1(M)[\bn^T].$
We denote the Riemannian metric on $N_1(M)[\bn^T]$ by $(G_{ij}(p,\bxi))_{1\leqslant i,j\leqslant n-1}.$
We now arbitrarily choose (at least locally) a unit
spacelike normal vector field $\bn^S$ with $\bn^S(u)\in N_1(M)_p[\bn^T]$, where $p=\bX(u).$
We call $(\bn^T,\bn^S)$ an {\it adopted pair of normal vector fields\/} along $M.$
Clearly, the vectors
$\bn^T (u)\pm \bn^S(u)$ are null. 
We define a mapping
\[
\mathbb{NG}(\bn^T,\bn^S):U\lon \Lambda^*
\]
by $\mathbb{NG}(\bn^T,\bn^S)(u)=\bn^T(u)+\bn^S(u).$
We call it the {\it nullcone Gauss image} of $M=\bX(U)$ with respect to
$(\bn^T,\bn^S).$
With the identification of $M$ and $U$ through $\bX,$ we have the
linear mapping provided by the derivative of the nullcone Gauss image $\mathbb{NG}(\bn^T,\bn^S)$ at each point $p\in M$ as follows:
\[
d_p\mathbb{NG}(\bn^T,\bn^S):T_pM\lon T_p\R^{n+1}_1= T_pM\oplus N_p(M).
\]
Consider the orthogonal projections $\pi ^t:T_pM\oplus
N_p(M)\rightarrow T_p(M).$ We define
\[
d_p\mathbb{NG}(\bn^T,\bn^S)^t=\pi ^t\circ d_p(\bn^T+\bn^S).
\]
We call the
linear transformation $S_{p}(\bn^T,\bn^S)=- d_p\mathbb{NG}(\bn^T,\bn^S)^t$ the {\it $(\bn^T,\bn^S)$-shape
operator} of $M=\bX (U)$ at $p=\bX (u).$ 
Let $\{\kappa
_{i}(\bn^T,\bn^S)(p)\}_ {i=1}^s$ be the eigenvalues of $S_{p}(\bn^T,\bn^S)$, which are called the {\it nullcone
principal curvatures  with respect to $(\bn^T,\bn^S) $\/} at $p=\bX(u)$.
Then the {\it nullcone Gauss-Kronecker curvature with respect to
$(\bn^T,\bn^S)$\/} at $p=\bX (u)$ is defined by
\[
K_N(\bn^T,\bn^S)(p)={\rm det} S_{p}(\bn^T,\bn^S).
\]
We say that a point $p=\bX (u)$ is an {\it $(\bn^T,\bn^S)$-umbilical
point} if 
\[
S_{p}(\bn^T,\bn^S)=\kappa (\bn^T,\bn^S)(p) 1_{T_{p}M}.
\]
We say that $M=\bX (U)$ is {\it totally
$(\bn^T,\bn^S)$-umbilical} if all points on $M$ are
$(\bn^T,\bn^S)$-umbilical.
Moreover, $M=\bX(U)$ is said to be {\it totally nullcone umbilical} if
it is totally $(\bn^T,\bn^S)$-umbilical for any adopted pair $(\bn^T,\bn^S).$
\par
We deduce now the nullcone Weingarten formula. Since $\bX _{u_i}$
$(i=1,\dots s)$ are spacelike vectors, we have a Riemannian metric
(the {\it first fundamental form \/}) on $M=\bX (U)$
defined by $ds^2 =\sum _{i=1}^{s} g_{ij}du_idu_j$,  where
$g_{ij}(u) =\langle \bX _{u_i}(u ),\bX _{u_j}(u)\rangle$ for any
$u\in U.$ We also have the {\it nullcone second fundamental invariant
with respect to the normal vector field $(\bn^T,\bn ^S) $\/} defined
by $h _{ij}(\bn^T,\bn^S )(u)=\langle -(\bn^T +\bn^S)
_{u_i}(u),\bX_{u_j}(u)\rangle$ for any $u\in U.$
By similar arguments to those in the proof of \cite[Proposition 3.2]{IzuSM}, we have 
the following proposition.
\begin{Pro}
We choose a pseudo-orthonormal frame $\{\bX, \bn^T,\bn^S_1,\dots ,\bn^S_{k-1}\}$ of $N(M)$ with $\bn^S_{k-1}=\bn^S.$ Then we have the following nullcone Weingarten formula {\rm :}
\vskip1.5pt
\par\noindent
{\rm (a)} $\mathbb{NG}(\bn^T,\bn^S)_{u_i}=\langle \bn ^T_{u_i},\bn ^S\rangle(\bn^T+\bn^S)+\sum _{\ell =1}^{k-2}\langle (\bn^T+\bn^S)_{u_i},\bn^S_\ell \rangle\bn^S_\ell -\sum_{j=1}^{s}
h_i^j(\bn^T,\bn^S )\bX _{u_j}$
\par\noindent
{\rm (b)} $
\pi ^t\circ \mathbb{NG}(\bn^T,\bn^S)_{u_i}=-\sum_{j=1}^{s}
h_i^j(\bn^T,\bn^S )\bX _{u_j}.
$
\smallskip
\par\noindent
Here $\displaystyle{\left(h_i^j(\bn^T,\bn^S )\right)=\left(h_{ik}(\bn^T,\bn^S)\right)\left(g^{kj}\right)}$
and $\displaystyle{\left( g^{kj}\right)=\left(g_{kj}\right)^{-1}}.$
\end{Pro}
\par
As a consequence of the above proposition, we have an explicit
expression of the nullcone curvature
by
$$
K_N (\bn^T,\bn^S )=\frac{\displaystyle{{\rm det}\left(h_{ij}(\bn^T,\bn^S )\right)}}
{\displaystyle{{\rm det}\left(g_{\alpha \beta}\right)}}.
$$
Since $\langle -(\bn^T +\bn^S )(u),\bX _{u_j}(u)\rangle =0,$ we have
$h_{ij}(\bn ^T,\bn^S)(u)=\langle \bn^T (u)+\bn^S (u),\bX
_{u_iu_j}(u)\rangle.$ Therefore the nullcone second fundamental
invariant at a point $p_0=\bX (u_0)$ depends only on the values 
$\bn^T (u_0)+\bn^S (u_0)$ and $\bX _{u_iu_j}(u_0)$, respectively.
Thus, the nullcone curvatures also depend only on
$\bn^T (u_0)+\bn^S (u_0)$, $\bX_{u_i}(u_0)$  and $\bX
_{u_iu_j}(u_0)$, which are independent of the derivations of the vector fields 
$\bn^T$ and $\bn^S .$ We write $\kappa _i(\bn^T_0,\bn^S_0)(p_0)$ $(i=1,\dots ,s)$
and $K_N (\bn
^T_0,\bn^S_0)(u_0)$ as the nullcone curvatures at $p_0=\bX (u_0)$
with respect to $(\bn ^T_0,\bn^S_0)=(\bn^T (u_0),\bn^S(u_0)).$ We
might also say that a point $p_0=\bX (u_0)$ is 
{\it $(\bn^T_0,\bn^S_0)$-umbilical\/} because the nullcone $(\bn^T,\bn^S) $-shape
operator at $p_0$ depends only on the normal vectors $(\bn
^T_0,\bn^S_0).$
So we denote that $h_{ij}(\bn^T,\bxi)(u_0)=h_{ij}(\bn^T,\bn^S)(u_0)$ and 
 $K_N(\bn^T,\bxi)(p_0)=K_N(\bn^T_0,\bn^S_0)(p_0)$,
where $\bxi =\bn^S(u_0)$ for some local extension $\bn^T(u)$ of $\bxi.$
Analogously, we say that a point $p_0=\bX (u_0)$ is an {\it $(\bn
^T_0,\bn^S_0)$-parabolic point \/} of $\bX :U\lon \R^{n+1}_1$ if
$K_N (\bn ^T_0,\bn^S_0)(u_0)=0.$ We also say that a point $p_0=\bX
(u_0)$ is a {\it $(\bn ^T_0,\bn^S_0)$-flat point \/} if it is an
$(\bn ^T_0,\bn^S_0)$-umbilical point and $K_N(\bn^T
_0,\bn^S_0)(u_0)=0.$
\par
On the other hand, we define a map
$
\mathbb{NG}(\bn^T):N_1(M)[\bn^T]\lon \Lambda ^*
$
by
$\mathbb{NG}(\bn^T)(u,\bxi)=\bn^T(u)+\bxi,$
which we call the {\it nullcone Gauss image} of $N_1(M)[\bn^T].$
This map leads us to the notions of curvatures.
Let $T_{(p,\xi)}N_1(M)[\bn^T]$ be the tangent space of $N_1(M)[\bn^T]$ at $(p,\bxi).$
Under the canonical identification $(\mathbb{NG}(\bn^T)^*T\R^{n+2}_2)_{(p,\sbxi)}
=T_{(\sbn^T(p)+\sbxi)}\R^{n+2}_2\equiv T_p\R^{n+2}_2,$
we have
\[
T_{(p,\sbxi)}N_1(M)[\bn^T]=T_pM\oplus T_\xi S^{k-2}\subset T_pM\oplus N_p(M)=T_p\R^{n+2}_2,
\]  
where $T_\xi S^{k-2}\subset T_\xi N_p(M)\equiv N_p(M)$ and $p=\bX(u).$
Let 
\[
\Pi ^t :\mathbb{NG}(\bn^T)^*T\R^{n+2}_2=TN_1(M)[\bn^T]\oplus \R^{k+1}
\lon TN_1(M)[\bn^T]
\]
be the canonical projection.
Then
we have a linear transformation
\[
S_N (\bn^T)_{(p,\sbxi)}=-\Pi^t_{\mathbb{LG}(n^T)(p,\xi)}\circ d_{(p,\xi)}\mathbb{NG}(\bn^T)
: T_{(p,\xi)}N_1(M)[\bn^T]\lon T_{(p,\xi)}N_1(M)[\bn^T],
\]
which is called the {\it nullcone shape operator} of $N_1(M)[\bn^T]$ at $(p,\bxi).$ 
Let $\kappa _N(\bn^T)_i(p,\bxi)$ be the eigenvalues of $S_N(\bn^T) _{(p,\sbxi)}$, $(i=1,\dots ,n-1)$. 
Here, we denote $\kappa _N(\bn^T)_i(p,\bxi)$, $(i=1,\dots ,s)$ as the eigenvalues belonging to
the eigenvectors on $T_pM$
and $\kappa _N(\bn^T)_i(p,\bxi)$, $(i=s+1,\dots n-1)$ as the eigenvalues belonging to the eigenvectors on 
the tangent space of the fiber  
of $N_1(M)[\bn^T].$  
Then we have the following proposition.
\begin{Pro}
We choose a {\rm (}local\/{\rm )} pseudo-orthonormal frame $\{\bX, \bn^T,\bn^S_1,\dots ,\bn^S_{k-1}\}$ of $N(M)$ with $\bn^S_{k-1}=\bn^S.$
For $p_0=\bX(u_0)$ and $\bxi_0=\bn^S(u_0),$ we have
$\kappa _N(\bn^T)_i(p_0,\bxi_0)=\kappa _i(\bn^T,\bn^S)(u_0)$, $(i=1,\dots ,s)$ and
$\kappa _N(\bn^T)_i(p_0,\bxi_0)=-1$, $(i=s+1,\dots n-1).$
\end{Pro}
\demo
Since $\{\bX, \bn^T,\bn^S_1,\dots ,\bn^S_{k-1}\}$ is a pseudo-orthonormal frame of$N(M),$ we have $\langle \bX(u_0),\bxi_0\rangle =
\langle \bn^T(u_0),\bxi_0\rangle =\langle \bn^S_i(u_0),\bxi_0\rangle =0.$
Therefore, we have $T_{\sbxi}S^{k-2}=\langle \bn^S_1(u_0),\dots ,\bn^S_{k-2}(u_0)\rangle .$
By this orthonormal basis of $T_{\sbxi_0}S^{k-2},$
the canonical Riemannian metric $G_{ij}(p_0,\bxi_0)$ is represented by
\[
(G_{ij}(p_0,\bxi))=\left(
\begin{array}{cc}
g_{ij}(p_0)  & 0 \\
0 & I_{k-2}
\end{array}
\right) ,
\]
where $g_{ij}(p_0)=\langle \bX_{u_i}(u_0), \bX_{u_j}(u_0)\rangle $.
\par
On the other hand, by Proposition 3.1, we have
\[
-\sum_{j=1}^s h^j_i(\bn^T,\bn^S)(u_0)\bX_{u_j}=\mathbb{NG}(\bn^T,\bn^S)_{u_i}(u_0)=
d_{p_0}\mathbb{NG}(\bn^T,\bn^S)\left(\frac{\partial}{\partial u_i}\right),
\]
so that we have
\[
S_N(\bn^T)_{(p_0,\bxi_0)}\left(\frac{\partial}{\partial u_i}\right)=\sum_{j=1}^s h^j_i(\bn^T,\bn^S)(u_0)\bX_{u_j}.
\]
Therefore, the representation matrix of $S_N(\bn^T)_{(p_0,\bxi_0)}$ with respect to the basis
$$
\{\bX_{u_1}(u_0),\dots ,\bX_{u_s}(u_0),\bn^S_1(u_0),\dots ,\bn^S_{k-2}(u_0)\}
$$ of $T_{(p_0,\bxi_0)}(N_1(M)[\bn^T])$
is of the form
\[
\left(
\begin{array}{cc}
h^j_i(\bn^T,\bn^S)(u_0)  & * \\
0 & -I_{k-2}
\end{array}
\right).
\]
Thus, the eigenvalues of this matrix are $\lambda _i=\kappa _i(\bn^T,\bn^S)(u_0)$, $(i=1,\dots ,s)$ and
$\lambda _i=-1$ , $(i=s+1,\dots ,n-1)$.
This completes the proof.
\enD
We call $\kappa _N(\bn^T)_i(p,\bxi)$, $(i=1,\dots ,s)$
the {\it nullcone principal curvatures} of $M$ with respect to $(\bn^T,\bxi)$ at $p\in M.$
The {\it nullcone Lipschitz-Killing curvature} of $N_1(M)[\bn^T]$ at $(p,\bxi)$
is defined to be $K_N (\bn^T)(p,\bxi)=\det S_N (\bn^T)_{(p,\xi)}.$

\section{Lightlike hypersurfaces in anti-de Sitter space}
 We define a hypersurface
$$
\mathbb{LH}_M(\bn^T):N_1(M)[\bn^T]\times \R\lon AdS^{n+1}
$$
by
$$
\mathbb{LH}_M(\bn^T)((p,\bxi),\mu)=\bX(u)+\mu (\bn^T+\bxi)(u)=\bX(u)+\mu\mathbb{NG}(\bn^T)(u,\bxi),
$$
where $p=\bX (u),$ which is called the {\it lightlike hypersurface\/} along $M$ relative to $\bn^T.$
In general, a hypersurface $H\subset AdS^{n+1}$ is called a {\it lightlike hypersurface\/} if it is tangent to
the lightcone at any regular point.
We remark that $\mathbb{NH}_M(\bn^T)(N_1(M)[\bn^T]\times\R)$ is a lightlike hypersurface.
 \par
 We introduce the notion of height
functions on spacelike submanifold, which is useful for the study of
singularities of lightlike hypersurfaces.
We define a family of functions 
$$H: M\times AdS^{n+1}\lon \R$$
 on a spacelike submanifold  $M=\bX (U)$ 
 by
$
 H(p,\blambda)=H(u,\blambda )=\langle \bX (u) ,\blambda\rangle +1,
$
where $p=\bX(u).$
We call
$H$ the {\it anti-de Sitter height function\/} (briefly, {\it AdS-height function\/}) on the spacelike submanifold
$M.$
For any fixed $\blambda _0\in AdS^{n+1},$ we write $h_{\sblambda _0}(p)=H(p,\blambda _0)$
and have
the following proposition.
\par
\begin{Pro}
Let $M$ be a spacelike submanifold 
and
$H: M\times(AdS^{n+1}\setminus M)\to\R$
the AdS-height function on $M.$
Suppose that $p_0=\bX(u_0)\not=\blambda _0.$ Then we have the following\/$:$
\par
{\rm (1)}
$h_{\sblambda _0}(p_0)=\partial h_{\sblambda _0}/\partial u_i(p_0)=0$, $(i=1,\dots ,s)$
if and only if
there exist $\bxi_0 \in N_1(M)_{p_0}[\bn^T]$ and $\mu_0\in
\R\setminus \{0\}$ such that 
$$
\blambda _0 =\bX(u_0)+\mu_0\mathbb{NG}(\bn^T)(u_0,\bxi_0)=\mathbb{LH}_M(\bn^T)((p_0,\bxi_0),\mu_0).
$$ 
\par
{\rm (2)}
$h_{\sblambda _0}(p_0)=\partial h_{\sblambda _0}/\partial u_i(p_0)=
{\rm det}{\mathcal H}(h_{\sblambda _0})(p_0)=0$ $(i=1,\dots ,s)$
if and only if
there exist $\bxi_0 \in N_1(M)_{p_0}[\bn^T]$ and $\mu_0\in
\R\setminus \{0\}$ such that 
$$
\blambda _0=\mathbb{LH}_M(\bn^T)((p_0,\bxi_0),\mu_0)
$$
and 
$1/\mu$ is one of the non-zero nullcone
principal curvatures 
$\kappa_N(\bn^T)_i(p_0,\bxi_0), (i=1,\dots ,s).$
\par
Here, ${\mathcal H}(h_{\sblambda _0})(p_0)$
is the Hessian matrix of $h_{\sblambda _0}$ at $p_0.$
\par
{\rm (3)} With condition {\rm (2)}, ${\rm rank}\, {\mathcal H}(h_{\sblambda _0})(p_0)=0$ if and only if
$p_0=\bX(u_0)$ is a non-flat $(\bn^T(u_0),\bxi _0)$-umbilical point.
\end{Pro}
\demo
(1) For $p=\bX(u),$ the condition $h_{\sblambda _0}(p)=\langle \bX(u),{\blambda
_0}\rangle+1 =0$
means that
\[
\langle \bX(u)-\lambda _0,\bX(u)-\blambda_0\rangle =\langle\bX(u),\bX(u)\rangle-2\langle\bX(u),\blambda _0\rangle+\langle\blambda _0,\blambda_0\rangle=-2(1+\langle\bX(u),\blambda _0\rangle)=0,
\]
so that
$\bX(u)-{\blambda _0}\in \Lambda^*.$
Since $\partial h_{\sblambda _0}/\partial u_i(p)=\langle \bX_{u_i}(u), {\blambda _0}\rangle $
and $\langle \bX_{u_i},\bX\rangle=0,$
we have $\langle \bX_{u_i}(u),\blambda _0\rangle=-\langle \bX_{u_i}(u)-\blambda _0\rangle$.
Therefore, $\partial h_{\sblambda _0}/\partial u_i(p)=0$
if and only if
$\bX (u)-{\blambda _0}\in N_pM.$
On the other hand, the condition $h_{\sblambda _0} (p)=\langle \bX(u),\blambda _0\rangle+1=0$
implies that $\langle\bX(u),\bX(u)-\blambda_0\rangle =0$.
This means that $\bX(u)-\blambda _0\in T_pAdS^{n+1}.$
Hence
$h_{\sblambda _0}(p_0)=\partial h_{\sblambda _0}/\partial u_i((p_0)=0$ $(i=1,\dots, s)$
if and only if
$\bX(u_0)-{\blambda _0}\in N_{p_0}M\cap \Lambda^*\cap T_{p_0}AdS^{n+1}.$
Then we denote that
$\bv=\bX (u_0)-{\blambda _0}\in N_{p_0}M\cap \Lambda^*\cap T_{p_0}AdS^{n+1}.$
If $\langle \bn^T(u_0),\bv\rangle =0,$ then $\bn^T(u_0)$ belongs to
a lightlike hyperplane in the Lorentz space $T_{p_0}AdS^{n+1},$ so that $\bn^T(u_0)$ is lightlike or spacelike.
This contradiction to the fact that $\bn^T(u_0)$ is a timelike unit vector. Thus,
$\langle \bn^T(u_0),\bv\rangle \not=0.$ 
We set
\[
\bxi_0=\frac{-1}{\langle \bn^T(u_0),\bv\rangle}\bv -\bn^T(u_0).
\]
Then we have
\begin{eqnarray*}
\langle \bxi_0,\bxi_0\rangle &=& -2\frac{-1}{\langle \bn^T(u_0),\bv\rangle} \langle \bn^T(u_0),\bv\rangle-1=1 \\
\langle \bxi_0,\bn^T(u_0)\rangle &=& \frac{-1}{\langle \bn^T(u_0),\bv\rangle} \langle \bn^T(u_0),\bv\rangle+1=0.
\end{eqnarray*}
This means that $\bxi_0\in N_1(M)_{p_0}(M)[\bn^T].$
Since $-\bv=\langle \bn^T(u_0),\bv\rangle(\bn^T(u_0)+\bxi_0),$
we have 
${\blambda _0}=\bX(u_0)+\mu_0\mathbb{NG}(\bn^T)(p_0,\bxi_0)$, where
$p_0=\bX(u_0)$ and $\mu_0=\langle \bn^T(u_0),\bv\rangle.$
For the converse assertion,  suppose that $\blambda_0=\bX(u_0)+\mu_0\mathbb{NG}(\bn^T)(p_0,\bxi_0).$
Then $\blambda_0-\bX(u_0)\in N_{p_0}(M)\cap \Lambda ^*$ and
$\langle\blambda_0-\bX(u_0),\bX(u_0)\rangle=\langle \mu_0\mathbb{NG}(\bn^T)(p_0,\bxi_0),\bX(u_0)\rangle=0.$
Thus we have $\blambda_0-\bX(u_0)\in N_{p_0}(M)\cap \Lambda ^*\cap T_{p_0}AdS^{n+1}.$
By the previous arguments, these conditions are equivalent to the condition that
$h_{\sblambda _0}(p_0)=\partial h_{\sblambda _0}/\partial u_i((p_0)=0$ $(i=1,\dots, s)$.

\par
(2) By a straightforward calculation, we have
\[
\frac{\partial ^2 h_{\sblambda_0}}{\partial u_i\partial u_j}(u)
=\langle\bX _{u_iu_j},\blambda _0\rangle.
\]
Under the conditions ${\blambda _0}=\bX(u_0)+\mu_0(\bn^T(u_0)+\bxi_0)$,
we have
\[
\frac{\partial ^2 h_{\sblambda_0}}{\partial u_i\partial u_j}(u_0)
=\langle \bX _{u_iu_j}(u_0),\bX(u_0)\rangle +\mu_0\langle\bX_{u_iu_j}(u_0), (\bn^T(u_0)+\bxi_0)\rangle .
\]
Since $\langle \bX_{u_i},\bX\rangle =0,$ we have $\langle\bX_{u_iu_j},\bX\rangle=-\langle \bX_{u_i},\bX_{u_j}\rangle.$
Therefore, we have
\[
\left(\frac{\partial ^2 h_{\sblambda_0}}{\partial u_i\partial u_\ell}(u_0)\right)\left(g^{j\ell}(u_0)\right)
=\left(\mu_0 h^j_i(\bn^T,\bn^S)(u_0)-\delta ^j_i\right),
\]
where $\bn^S$ is the local section of $N_1(M)[\bn^T]$ with $\bn^S(u_0)=\bxi _0.$
It follows that ${\rm det}{\mathcal H}(g)(p_0)=0$ if and only if 
$1/\mu_0$ is an eigenvalue of $(h^i_j(\bn^T,\bn^S)(u_0)),$ which is equal to
one of the nullcone principal curvatures $\kappa _i(\bn^T,\bn^S)(u_0)=\kappa_N(\bn^T)_i(p_0,\bxi_0), (i=1,\dots ,s)$.
\par
(3) By the above calculation,  ${\rm rank}\, {\mathcal H}(h_{\sblambda _0})(p_0)=0$ if and only if
$$
(h^i_j(\bn^T,\bn^S)(u_0))=\frac{1}{\mu _0}(\delta ^j_i),
$$
where $1/\mu_0=\kappa _N(\bn^T)_i(p_0,\bxi_0),\ (i=1,\dots ,s)$. This means that $p_0=\bX(u_0)$ is an $(\bn^T(u_0),\bxi_0)$-umbilical point.
\enD
\par
In order to understand the geometric meaning of the assertions of Proposition 4.1, we briefly review the theory of Legendrian singularities.
For detailed expressions, see \cite{Arnold1, Zak}.
Let 
$\pi :PT^*(\R^{n+1}) \longrightarrow \R^{n+1}$ be the projective cotangent bundle with its canonical contact structure.
We  next review the geometric properties of this bundle.
Consider the tangent bundle
$
\tau :TPT^*(\R^{n+})\rightarrow PT^*(\R^{n+1})
$
and the differential map
$
d\pi :TPT^*(\R^{n+1})\rightarrow T\R^{n+1}
$
of $\pi .$
For any $X\in TPT^*(\R^{n+1}),$ there exists an element
$\alpha\in T^*(\R^{n+1}_1$ such that
$\tau (X)=[\alpha ].$  For an element $V\in T_x(\R^{n+1}),$
the property $\alpha (V)=0$ does not depend on the choice of
representative of the class $[\alpha ].$  Thus we can define the canonical
contact structure on $PT^*(\R^{n+1})$ by
$$
K=\{X\in TPT^*(\R^{n+1})\ |\ \tau (X)(d\pi (X))=0\}.
$$
We have a trivialization
$
PT^*(\R^{n+1})\cong
\R^{n+1}\times P^n(\R)^*,
$
and call
$$
((v_0,v_1,\dots ,v_n),[\xi _0:\xi _1:\cdots :\xi _n])
$$
the {\it  homogeneous coordinates} of $PT^*(\R^{n+1}),$ where
$
[\xi _0:\xi _1:\cdots :\xi _n]
$
are the homogeneous coordinates of the dual projective space
$P^n(\R)^*.$
It is easy to show that $X\in K_{(x,[\xi])}$ if and only if
$
\sum_{i=0}^n \mu _i\xi _i=0,
$
where
$
d\tilde\pi (X)=\sum_{i=0}^n \mu _i\partial/\partial v_i.
$
An immersion  $i:L\rightarrow PT^*(\R^{n+1})$ is said to be
{\it a Legendrian immersion} if $\text{dim}\, L=n$ and $di_q(T_qL)\subset
K_{i(q)}$
for any $q\in L.$
The map $\pi\circ i$ is also called {\it the Legendrian map} and the set
$W(i)=\text{image}\, \pi\circ i$, the {\it wave front set} of $i.$
Moreover, $i$ (or, the image of $i$) is called the {\it Legendrian lift }
of $W(i).$
\par
Let $F:({\mathbb R}^k\times{\mathbb R}^{n+1},\bo )\longrightarrow ({\mathbb
R},\bo )$ be a
function germ.
We say that $F$ is {\it a Morse family of hypersurfaces} if the map germ
$$
\Delta^*F=\left(F,\frac{\partial F}{\partial q_1},\dots ,\frac{\partial
F}{\partial q_k}
\right):({\mathbb R}^k\times {\mathbb R}^{n+1},\bo )\lon ({\mathbb R}\times
{\mathbb R}^k,\bo )
$$
is submersive, where $(q,x)=(q_1,\dots ,q_k,x_0,\dots ,x_n)\in  ({\mathbb
R}^k\times
{\mathbb R}^{n+1},\bo ).$
In this case we have a smooth $n$-dimensional submanifold
$$
\Sigma _*(F)=\Bigl\{(q,x)\in ({\mathbb R}^k\times{\mathbb R}^{n+1},\bo )\  \Bigm|
\ F(q,x)=\frac{\partial F}{\partial q_1}(q,x)=\cdots =\frac{\partial F}{\partial q_k}(q,x)=0
\ \Bigr\}
$$
and the map germ $\mathscr{L} _F:(\Sigma _*(F), \bo)\lon PT^*{\mathbb R}^{n+1}$ defined by
$$
\mathscr{L} _F(q,x)=\left(x,\left[\frac{\partial F}{\partial x_0}(q,x):\cdots :\frac{\partial F}{\partial x_n}(q,x)\right]\right)
$$
is a Legendrian immersion. 
We call $F$ {\it a generating family} of $\mathscr{L} _F(\Sigma _*(F)),$
and the wave front set is given by
$
W(\mathscr{L} _F)\! =\pi _n(\Sigma _*(F)),
$
where $\pi _n:\R^k\times\R^n\lon \R^n$ is the canonical projection.
In the theory of unfoldings of function germs, the wave front set $W(\mathscr{L} _F)$ is called a
{\it discriminant set} of $F,$ which we also denote $\mathcal{D}_F.$
Therefore, Proposition 4.1 asserts that the discriminant set of the
AdS-height function $H$ is given by
\[
{\mathcal D}_{H}=\Bigl\{\blambda\in AdS^{n+1}  \Bigm| \blambda =\bX (u)+\mu (\bn^T\pm \bxi)
(u),\ p=\bX(u)\in M, \bxi\in N_1(M)_p[\bn^T], \mu \in \R\ \Bigr\},
\]
which is the image of the lightlike hypersurface along $M$ relative to $\bn^T.$
\par
By the assertion (2) of Proposition 4.1, a singular point of the lightlike hypersurface is
a point $\blambda _0=\bX(u_0)+\mu _0(\bn^T+\bxi_0)(u_0)$
for
$p_0=\bX(u_0)$ and 
$\mu _0 =1/\kappa  _N(\bn^T)_i(p_0,\bxi_0),$ $i=1,\dots .s).$
Then we have the following corollary.
\begin{Co}
 The critical value of $\mathbb{LH}_M(\bn^T)$ is the point 
 \[
\blambda =\bX(u)+\frac{1}{\kappa _N(\bn^T)_i(p,\bxi)}\mathbb{LG}(\bn^T)(u,\bxi),
\]
where
$p=\bX(u)$ and 
$\kappa_N(\bn^T)_i(p,\bxi)\not= 0.$
\end{Co}
\par
For a non-zero nullcone principal curvature $\kappa _N(\bn^T)_i(p_0,\bxi_0)\not= 0,$ we
have an open subset $O_i\subset N_1(M)[\bn^T]$ such that $\kappa _N(\bn^T)(p,\bxi)\not= 0$
Therefore, we have a non-zero nullcone principal curvature function $\kappa _N(\bn^T):O_i\lon \R$.
We define a mapping
$
\mathbb{LF}_{\kappa_N(\bn^T)_i} :O_i\lon AdS^{n+1}
$
by
\[
\mathbb{LF}_{\kappa_N(\sbn^T)_i}(p,\bxi)=\bX(u)+\frac{1}{\kappa_N(\bn^T)_i(p,\bxi)}\mathbb{NG}(\bn^T)(u,\bxi),
\]
where $p=\bX(u).$
We also  define
\[
\mathbb{LF}_M(\bn^T)=\bigcup\left\{\mathbb{LF}_{\kappa_N(\sbn^T)_i}(p,\bxi)\ |\ (p,\bxi)\in N_1(M)[\bn^T]\ \mbox{s.t.}\ \kappa _N(\bn^T)_i(p,\bxi)\not= 0,i=1,\dots ,s \right\} .
\]
We call $\mathbb{LF}_{M}(\bn^T)$ the {\it lightlike focal set} of $M=\bX(U)$
relative to $\bn^T.$
By definition, the lightlike focal set of $M=\bX(U)$
relative to $\bn^T$ is the critical values set of the lightlike hypersurface
$\mathbb{LH}_M(\bn^T)(N_1(M)[\bn^T]\times\R)$ along $M$ relative to $\bn^T.$

\par
We can show that the image of the lightlike hypersurface along $M$ is independent of the choice of the future directed
timelike normal vector field $\bn^T$ as a corollary of Proposition 4.1.
\begin{Co}
Let $\bn^T$ and $\overline{\bn}^T$ be future directed timelike unit normal fields along $M$.
Then we have
\[
\mathbb{LH}_M(\bn^T)(N_1(M)[\bn^T]\times\R)=\mathbb{LH}_M(\overline{\bn}^T)(N_1(M)[\overline{\bn}^T]\times\R)\ \mbox{and}\ \mathbb{LF}_M(\bn^T)=\mathbb{LF}_M(\overline{\bn}^T).
\]
\end{Co}
\demo
By Proposition 4.1, the images of the lightlike hypersurface along $M$ relative to $\bn^T$ and $\overline{\bn}^T$ are
the discriminant sets of the  AdS-height function $H$ on $M$.
Moreover, the focal set is the critical value set of the lightlike hypersurface
along $M$ relative to $\bn^T.$
Since $H$ is independent of the choice of $\bn^T,$ we have the assertion.
\enD
\par
We have the following proposition.
\begin{Pro}
Let $H$ be the AdS-height function on $M.$
For any point $(u,\blambda )\in \Delta ^*H^{-1}(0),$ the germ of $H$ at
$(u,\blambda )$
is a Morse family of hypersurfaces. 
\end{Pro}
\demo
We denote that
$$
\bX(u)=(X_{-1}(u),X_0(u),X_1(u),\dots ,X_n(u))\ {\rm and}\
\blambda =(\lambda _{-1},\lambda _0,\lambda _1,\dots ,\lambda _n).
$$
We define an open subset 
$U_{-1}^+=\{\blambda \in AdS^{n+1}\ |\ \lambda _{-1}>0\ \}$.
For any $\blambda \in U_{-1}^+,$ we have
\[
\lambda _{-1}=\sqrt{1-\lambda _0^2+\lambda _1^2+\cdots \lambda _n^2}.
\]
Thus, we have a local coordinate of $AdS^{n+1}$ given by $(\lambda _0,\lambda _1,\dots ,\lambda _n)$ on
$U_{-1}^+.$
By definition, we have
$$
H(u,\blambda )=-X_{-1}(u)\sqrt{1-\lambda _0^2+\sum_{i=1}^n \lambda _i^2}-X_0(u)\lambda_0+
X_1(u)\lambda_1 +\cdots 
+X_n(u)\lambda _n.
$$
We now prove that the mapping $$
\Delta^*H=\left(H, \frac{\partial H}{\partial u_1},\dots ,\frac{\partial H}{\partial
u_s}\right)
$$
is non-singular at $(u,\blambda )\in \Delta ^*H^{-1}(0).$
Indeed, the Jacobian matrix of $\Delta ^*H$ is given by
\newfont{\bg}{cmr10 scaled\magstep5}
\newcommand{\bigA}{\smash{\lower1.0ex\hbox{\bg A}}}
\[
\left(
\begin{array}{ccccc}
 &
X_{-1}\displaystyle{\frac{\lambda _0}{\lambda _{-1}}}-X_0 & -X_{-1}\displaystyle{\frac{\lambda _1}{\lambda _{-1}}}+X_1 &  \cdots  &
-X_{-1}\displaystyle{\frac{\lambda _n}{\lambda _{-1}}}-X_n \\
\bigA & X_{-1u_1}\displaystyle{\frac{\lambda _0}{\lambda _{-1}}}-X_{0u_1} &
-X_{-1u_1}\displaystyle{\frac{\lambda _1}{\lambda _{-1}}}+X_{1u_1}& \cdots  &-X_{-1u_1}\displaystyle{\frac{\lambda _n}{\lambda _{-1}}}-X_{nu_1}\\
&\vdots & \vdots & \ddots & \vdots \\
& X_{-1u_s}\displaystyle{\frac{\lambda _0}{\lambda _{-1}}}-X_{0u_s} &
-X_{-1u_s}\displaystyle{\frac{\lambda _1}{\lambda _{-1}}}+X_{1u_s}& \cdots  &-X_{-1u_s}\displaystyle{\frac{\lambda _n}{\lambda _{-1}}}-X_{nu_s}
\end{array}
\right) ,
\]
where
\begin{eqnarray*}
\bigA=
\left(\!\!
\begin{array}{ccc}
\langle \bX_{u_1} ,\blambda\rangle & \!\! \cdots\!\! & \langle \bX_{u_s},\blambda \rangle \\
\langle \bX_{u_1u_1},\blambda\rangle & \!\! \cdots\!\!  &
\langle \bX_{u_1u_s},\blambda \rangle \\
\vdots & \!\!\ddots\!\! & \vdots \\
\langle \bX_{u_su_1},\blambda \rangle &\!\! \cdots\!\! &
\langle \bX_{u_su_s},\blambda\rangle 
\end{array}
\!\!
\right) .
\end{eqnarray*}
\newcommand{\bigB}{\smash{\lower1.0ex\hbox{\bg B}}}
We now show that
the rank of 
\[
\bigB=
\left(
\begin{array}{cccc}
X_{-1}\displaystyle{\frac{\lambda _0}{\lambda _{-1}}}-X_0 & -X_{-1}\displaystyle{\frac{\lambda _1}{\lambda _{-1}}}+X_1 &  \cdots  &
-X_{-1}\displaystyle{\frac{\lambda _n}{\lambda _{-1}}}-X_n \\
 X_{-1u_1}\displaystyle{\frac{\lambda _0}{\lambda _{-1}}}-X_{0u_1} &
-X_{-1u_1}\displaystyle{\frac{\lambda _1}{\lambda _{-1}}}+X_{1u_1}& \cdots  &-X_{-1u_1}\displaystyle{\frac{\lambda _n}{\lambda _{-1}}}-X_{nu_1}\\
\vdots & \vdots & \ddots & \vdots \\
 X_{-1u_s}\displaystyle{\frac{\lambda _0}{\lambda _{-1}}}-X_{0u_s} &
-X_{-1u_s}\displaystyle{\frac{\lambda _1}{\lambda _{-1}}}+X_{1u_s}& \cdots  &-X_{-1u_s}\displaystyle{\frac{\lambda _n}{\lambda _{-1}}}-X_{nu_s}
\end{array}
\right) 
\]
is $s+1$ at $(u,\blambda)\in \Sigma _*(H).$
Since $(u,\blambda)\in \Sigma _*(H),$ we have
\[
\blambda =\bX(u)+\mu\left(\bn^T(u)+\sum _{i=1}^{k-1}\xi _i\bn_i(u)\right)
\]
with $\sum_{i=1}^{k-1}\xi ^2_i=1,$ where 
$\{\bX,\bn^T,\bn^S_1,\dots ,\bn^S_{k-1}\}$ is a pseudo-orthonormal (local) frame of $N(M).$
Without the loss of generality, we assume that $\mu\not= 0$ and $\xi_{k-1}\not= 0.$
We denote that
\[
\bn^T(u)=^t\!\!(n^T_{-1}(u),n^T_0(u),\dots n^T_n(u)),\ 
\bn_i(u)=^t\!\!(n^i_{-1}(u),n_0^i(u),\dots n^i_n(u)).
\]
\newcommand{\bigC}{\smash{\lower1.0ex\hbox{\bg C}}}
It is enough to show that the rank of 
the matrix
\[
\bigC =
\left(
\begin{array}{cccc}
X_{-1}\displaystyle{\frac{\lambda _0}{\lambda _{-1}}}-X_0 & -X_{-1}\displaystyle{\frac{\lambda _1}{\lambda _{-1}}}+X_1 &  \cdots  &
-X_{-1}\displaystyle{\frac{\lambda _n}{\lambda _{-1}}}-X_n \\
 X_{-1u_1}\displaystyle{\frac{\lambda _0}{\lambda _{-1}}}-X_{0u_1} &
-X_{-1u_1}\displaystyle{\frac{\lambda _1}{\lambda _{-1}}}+X_{1u_1}& \cdots  &-X_{-1u_1}\displaystyle{\frac{\lambda _n}{\lambda _{-1}}}-X_{nu_1}\\
\vdots & \vdots & \ddots & \vdots \\
 X_{-1u_s}\displaystyle{\frac{\lambda _0}{\lambda _{-1}}}-X_{0u_s} &
-X_{-1u_s}\displaystyle{\frac{\lambda _1}{\lambda _{-1}}}+X_{1u_s}& \cdots  &-X_{-1u_s}\displaystyle{\frac{\lambda _n}{\lambda _{-1}}}-X_{nu_s} \\
n^T_{-1}\displaystyle{\frac{\lambda _0}{\lambda _{-1}}}-n^T_0 & -n^T_{-1}\displaystyle{\frac{\lambda _1}{\lambda _{-1}}}+n^T_1 &  \cdots  &
-n^T_{-1}\displaystyle{\frac{\lambda _n}{\lambda _{-1}}}-n^T_n \\
n^1_{-1}\displaystyle{\frac{\lambda _0}{\lambda _{-1}}}-n^1_0 & -n^1_{-1}\displaystyle{\frac{\lambda _1}{\lambda _{-1}}}+n^1_1 &  \cdots  &
-n^1_{-1}\displaystyle{\frac{\lambda _n}{\lambda _{-1}}}-n^1_n \\
\vdots & \vdots & \ddots & \vdots \\
n^{k-2}_{-1}\displaystyle{\frac{\lambda _0}{\lambda _{-1}}}-n^{k-2}_0 & -n^{k-2}_{-1}\displaystyle{\frac{\lambda _1}{\lambda _{-1}}}+n^{k-2}_1 &  \cdots  &
-n^{k-2}_{-1}\displaystyle{\frac{\lambda _n}{\lambda _{-1}}}-n^{k-2}_n 
\end{array}
\right) 
\]
is $n+1$ at $(u,\blambda)\in \Sigma _*(H).$
We denote that
\[
\ba_i=^t\!\!(x_i(u),x_{iu_1}(u),\dots x_{iu_s}(u),n^T_i(u),n^1_i(u),\dots ,n^{k-2}_i(u)).
\]
Then we have
\[
\bigC=\left(\ba_{-1}\frac{\lambda_0}{\lambda _{-1}}-\ba_0,-\ba_{-1}\frac{\lambda _1}{\lambda _{-1}}+\ba_1,\dots ,
-\ba _{-1}\frac{\lambda _n}{\lambda _{-1}}+\ba_n\right).
\]
It follows that
\begin{eqnarray*}
\det \bigC\!\!\!\!\!\!\!
&{}&=\frac{\lambda _{-1}}{\lambda _{-1}}\det (\ba_0,\ba_1,\dots,\ba_n)+\frac{\lambda_0}{\lambda_{-1}}\det (\ba_{-1}\ba_{1},\dots ,\ba_n)
\\
&{}&-\frac{\lambda _1}{\lambda _{-1}}(-1)\det(\ba_{-1},\ba_0,\ba_2,\dots ,\ba_n)-\cdots 
-\frac{\lambda_n}{\lambda _{-1}}(-1)^{n-1}\det(\ba _{-1}\ba_0,\ba_1,\dots ,\ba_{n-1}).
\end{eqnarray*}
Moreover, we define
$\delta _i=\det (\ba_{-1},\ba_0,\ba_1,\dots,\ba _{i-1},\ba_{i+1},\dots ,\ba_n)$ for $i=-1,0,1,\dots ,n$
and
$\ba=(-\delta _{-1},-\delta _0,-\delta _1,(-1)^2\delta _2,\dots ,(-1)^{n-1}\delta _n).$
Then we have
\[
\ba=\bX\wedge\bX_{u_1}\wedge\cdots \wedge \bX_{u_s}\wedge\bn^T\wedge\bn_1\wedge\cdots\wedge\bn_{k-2}.
\]
We remark that $\ba\not=0$ and $\ba=\pm\|\ba\|\bn_{k-1}.$
By the above calculation, we have
\begin{eqnarray*}
\det\bigC\!\!\!\!\!\!\!&{}&=\left\langle\left(\frac{\lambda{-1}}{\lambda_{-1}},\frac{\lambda_0}{\lambda _{-1}},\dots,\frac{\lambda_n}{\lambda _{-1}}\right), \ba\right\rangle =\frac{1}{\lambda_{-1}}\left\langle \bX(u)+\mu\left(\bn^T(u)+\sum_{i=1}^{k-1}\xi_i\bn_i(u)\right),\ba\right\rangle \\
&{}&=\frac{1}{\lambda _{-1}}\times\pm\mu\xi_{k-1}\|\ba\|=\pm \frac{\mu\xi_{k-1}\|\ba\|}{\lambda _{-1}}\not= 0.
\end{eqnarray*}
Therefore the Jacobi matrix of $\Delta^*H$
is non-singular at $(u,\blambda )\in \Delta ^*H^{-1}(0).$
\par
For other local coordinates of $AdS^{n+1}$, we can apply the same method for the proof as the above case.
This completes the proof.
\enD
\par
Here we also consider the local coordinate 
$U^+_{-1}$.
Since $H$ is a Morse family of hypersurfaces, we have a Legendrian immersion
\[
\mathscr{L} _H:\Sigma _*(H)\lon PT^*(AdS^{n+1})|U^+_{-1}
\]
by the general theory of Legendrian singularities.
By definition, we have
\[
\frac{\partial H}{\partial \lambda _0}(u,\blambda) \\
=X_{-1}(u)\displaystyle{\frac{\lambda _0}{\lambda _{-1}}}-X_0(u),\ \frac{\partial H}{\partial \lambda _i}(u,\blambda)=-X_{-1}(u)\displaystyle{\frac{\lambda _i}{\lambda _{-1}}}+X_i(u),\ (i=1,\dots ,n).
\]
It follows that
\begin{eqnarray*}
&{}&\left[\frac{\partial H}{\partial \lambda _0}(u,\blambda):\frac{\partial H}{\partial \lambda _1}(u,\blambda):\cdots :\frac{\partial H}{\partial \lambda _n}(u,\blambda)\right]
\\
&{}&\qquad =[X_{-1}(u)\lambda _0-X_0(u)\lambda _{-1}:X_1(u)\lambda _{-1}-X_{-1}(u)\lambda _1:\cdots :X_n(u)\lambda _{-1}-X_{-1}(u)\lambda _n].
\end{eqnarray*}
Therefore, we have
\[
\mathscr{L}_H(u,\blambda )=(\blambda, [X_{-1}(u)\lambda _0-X_0(u)\lambda _{-1}:X_1(u)\lambda _{-1}-X_{-1}(u)\lambda _1:\cdots :X_n(u)\lambda _{-1}-X_{-1}(u)\lambda _n]),
\]
where
\[
\Sigma _*(H)=\{(u,\blambda)\ |\ \blambda =\mathbb{LH} _M(\bn^T)(p,\bxi,t)\ ((p,\bxi),t)\in N_1(M)[\bn^T]\times\R\}.
\]
We observe that $H$ is a generating family of the Legendrian submanifold $\mathscr{L}_H(\Sigma _*(H))$ whose wave front is $\mathbb{LH} _M(\bn^T)(N_1(M)[\bn^T]\times\R)$.
Therefore we say that the AdS-height 
function
$H$ on $M$ gives an {\it AdS-canonical generating family for the Legendrian lift\/} of $\mathbb{LH} _M(\bn^T)(N_1(M)[\bn^T]\times\R)$.
For other local coordinates of $AdS^{n+1},$ we have the similar results to the above case.

\section{Contact with lightcones}
\par
In this section we consider the geometric meaning of the singularities of lightlike hypersurfaces in Anti-de Sitter space from the view point of the theory of contact of submanifolds with model hypersurfaces in \cite{Montaldi}.
We begin with the following basic observations.
\begin{Pro}
Let $\blambda _0\in AdS^{n+1}$ and $M=\bX(U)$ a spacelike submanifold without points
satisfying $K_N (\bn^T)(p,\bxi)= 0.$
Then $M\subset
\Lambda ^{n+1}_{\lambda _0}\cap AdS^{n+1}$ 
if and only if $\{\blambda _0\}=\mathbb{LF}_M(\bn^T)$ is the lightcone focal set.
In this case we have $\mathbb{LH}_M(\bn^T)(N_1(M)[\bn^T])\subset \Lambda ^{n+1}_{\lambda _0}\cap AdS^{n+1}$
and $M=\bX(U)$ is totally nullcone umbilical.
\end{Pro}
\demo
By Proposition 3.1, $K_N(\bn^T)(p_0,\bxi_0)\not= 0$ if and only if
\[
\{(\bn^T+\bn^S), (\bn^T+\bn^S)_{u_1},
\dots , (\bn^T+\bn^S)_{u_{s}}\}
\]
is linearly independent for  $p_0=\bX(u_0)\in M$ and $\bxi_0=\bn^S(u_0),$
where $\bn^S:U\lon N_1(M)[\bn^T]$ is a local section.
By the proof of the assertion (1) of Proposition 4.1, 
$M\subset \Lambda ^{n+1}_{\lambda _0}\cap AdS^{n+1}$
if and only if $h_{\sblambda _0}(u)= 0$ for any $u\in U,$
where $h_{\sblambda _0}(u)=H(u,\blambda _0)$ is the AdS-height function on $M.$
It also follows from  Proposition 4.1 that there exists a smooth function $\eta 
:U\times N_1(M)[\bn^T]\lon \R$ and section $\bn^S:U\lon N_1(M)[\bn^T]$ such that
\[
\bX(u)=\blambda _0+\eta (u,\bn^S(u))(\bn^T(u)\pm\bn^S(u)).
\]
In fact, we have $\eta (u,\bn^S(u))=-1/\kappa _N(\bn^T)_i(p,\bxi)$ $i=1,\dots, s$, where
$p=\bX(u)$ and $\bxi=\bn^S(u).$
It follows that $\kappa _N(\bn^T)_i(p,\bxi)=\kappa _N(\bn^T)_j(p,\bxi),$ so that
$M=\bX(U)$ is totally nullcone umbilical.
Therefore we have
\[
\mathbb{LH}_M(\bn^T)(u,\bn^S(u),\mu)=\blambda _0+(\mu +\eta (u,\bn^S(u))(\bn^T(u)\pm\bn^S(u)).
\]
Hence we have
$\mathbb{LH}_M(\bn^T) (N_1(M)[\bn^T]\times \R)\subset \Lambda ^{n+1}_{\lambda _0}.$
By Corollary 4.2, the critical value set of $\mathbb{LH}_M(\bn^T) (N_1(M)[\bn^T]\times \R)$ is the lightlike focal set
$\mathbb{LF}_M(\bn^T).$
However, it is equal to $\lambda _0$ by the previous arguments. 
\par
For the converse assertion, suppose that $\blambda _0=\mathbb{LF}_M(\bn^T).$
Then we have 
\[
\blambda _0=\bX(u)+\frac{1}{\kappa _N(\bn^T)_i(\bX(u),\bxi)}\mathbb{LG}(\bn^T)(u,\bxi),
\]
for any $i=1,\dots ,s$ and $(p,\bxi)\in N_1(M)[\bn^T],$ where $p=\bX(u).$
Thus, we have 
\[
\kappa _N(\bn^T)_i(\bX(u),\bxi)=\kappa _N(\bn^T)_j(\bX(u),\bxi)
\]
for any $i,j=1,\dots ,s,$ so that $M$ is totally nullcone umbilical.
Since $\mathbb{LG}(\bn^T)(u,\bxi)$ is null, we have $\bX(u)\in \Lambda ^{n+1}_{\lambda_0}.$
This completes the proof.
\enD
\par
According to the above proposition, $\Lambda ^{n+1}_{\sblambda _0}\cap AdS^{n+1}$ is regarded as a model lightlike hypersurface in $AdS^{n+1}.$
We define
\[
T(AdS^{n+1})_{\sblambda _0}=\{\bx\in \R^{n+2}_2\ |\ \ \bx-\blambda _0\in T_{\sblambda_0}AdS^{n+1}\ \},
\]
where $T_{\sblambda_0}AdS^{n+1}$ is the tangent space of $AdS^{n+1}$ at $\blambda _0\in AdS^{n+1}.$
We call $T(AdS^{n+1})_{\sblambda _0}$ a {\it tangent affine space} of $AdS^{n+1}$ at $\blambda _0\in AdS^{n+1}.$
It is easy to show that
\[
\Lambda ^{n+1}_{\sblambda _0}\cap AdS^{n+1}=T(AdS^{n+1})_{\sblambda _0}\cap AdS^{n+1}.
\]
We denote that $LC_{\sblambda _0}(AdS^{n+1})=\Lambda ^{n+1}_{\sblambda _0}\cap AdS^{n+1}=T(AdS^{n+1})_{\sblambda _0}\cap AdS^{n+1}$,
which is called an {\it AdS-lightcone} with the vertex $\blambda _0\in AdS^{n+1}.$
Therefore, the model lightlike hypersurface is an AdS-lightcone.
\par
We consider the contact of
spacelike submanifolds with AdS-lightcones. 
Let 
\[
{\mathcal H}:AdS^{n+1}\times AdS^{n+1}\longrightarrow \R
\]
be a function defined
by ${\mathcal H}(\bx,\blambda)=\langle\bx , \blambda
\rangle+1 .$
Given $\blambda _0\in AdS^{n+1},$
we denote  ${\mathfrak h}_{\lambda _0}(\bx)={\mathcal H}(\bx ,\blambda _0)$,
so that we have ${\mathfrak h}_{\lambda _0}^{-1}(0)=
LC_{\sblambda _0}(AdS^{n+1}).$
For any $p_0=\bX(u_0)\in M$, $\mu_0\in \R$ and $\bxi_0\in N_1(M)_p[\bn^T],$ we consider the point
$\blambda _0=\bX(u_0)+\mu_0(\bn^T(u_0)+ \bxi_0).$
Then we have
$$
{\mathfrak h}_{\lambda _0}\circ\bX (u_0))={\mathcal H}\circ
(\bX\times 1_{AdS^{n+1}})(u_0,\blambda _0)
=H(p_0,\blambda _0)=0,
$$
where $\mu_0=1/\kappa_N(\bn^T)_i(p_0,\bxi_0),$ $i=1,\dots , s.$
We also have relations
$$
\frac{\partial {\mathfrak h}_{\lambda _0}\circ\bX}{\partial u_i}(u_0)=
\frac{\partial H}{\partial u_i}(p_0,\blambda _0)=0,\
i=1,\dots ,s.
$$
These imply that the AdS-lightcone ${\mathfrak h}_{\lambda_0}^{-1}(0)=
LC_{\sblambda_0}(AdS^{n+1})$ is tangent to $M=\bX (U)$  at  $p_0=\bX (u_0).$
In this case, we call $LC_{\sblambda _0}(AdS^{n+1})$ a {\it tangent AdS-lightcone} of
$M=\bX(U)$ at $p_0=\bX(u_0),$ which is denoted by $TLC_{\sblambda _0}(M)_{p_0}.$
Moreover, the tangent AdS-lightcone $TLC_{\sblambda _0}(M)_{p_0}$ is called an
{\it osculating AdS-lightcone} if $\blambda _0=\mathbb{LF}_{\kappa _N(\sbn^T)_i(p_0,\sbxi_0)}(u_0)\in \mathbb{LF}_M,$ for one of the nullcone principal curvature $\kappa _N(\bn^T)_i(p_0,\bxi_0).$
In this case, we call $\blambda _0$ the {\it center of the nullcone principal curvature} $\kappa _N(\bn^T)_i(p_0,\bxi_0)(u_0).$
Therefore, we can interpret that the lightlike focal set is the locus of the centers of 
nullcone principal curvatures. This fact is analogous to the notion of the focal sets of submanifolds in Euclidean space.

\par
Firstly, we consider a special contact of $M=\bX(U)$ with AdS-lightcones.
We say that 
$p_0=\bX(u_0)$ is  an {\it AdS-lightlike $k$-ridge  point} if $h_{\sblambda _0}$ 
has the $A_{k+2}$ singular point at $u_0$ for some $k\ge 1,$
where $\blambda _0 \in \mathcal{D}_{H}=\mathbb{LH}_M(N_1(M)\times \R).$ 
We simply say that  $p_0=\bX(u_0)$ is  an {\it AdS-lightlike ridge  point} if it is the AdS-lightlike $k$-ridge point 
for some $k\ge 1.$
For a function germ $f:(\R^{s},\widetilde{u} _0) \lon \R$, 
$f$ has an {\it $A_{k}$ singular point} at $\widetilde{u} _0$ if $f$ is 
$\mathcal{K}$-equivalent to the germ 
$ u^{k+1}_1 \pm u^2_2 \pm \dots \pm u^2_{s}$.
We say that two function germs $f_i:(\R^{s},\widetilde{u} _i)\lon \R$ $(i=1,2)$ are 
${\mathcal K}$-{\it equivalent} if there exists a diffeomorphism germ 
$\Phi :(\R^{s},\widetilde{u} _1)\lon  (\R^{s},\widetilde{u} _2)$ and  a non-zero function germ $\rho:(\R^{s},\widetilde{u}_1)\lon \R$ such that
$f_2\circ\Phi (u)=\rho (u)f_1(u).$ We consider the geometric meaning of AdS-lightlike ridge points.
Let $F: AdS^{n+1} \lon \R$ be a submersion and $\bX:U \lon AdS^{n+1}$ a spacelike embedding from an open set
$U\subset \R^s.$
We say that $M=\bX(U)$ and $F^{-1}(0)$ have a {\it corank $r$ contact} at $p_0=\bX(u_0)$ 
if the Hessian of the function $g(u) = F \circ \bX(u) $ has corank $r$ at $u_0$. 
We also say that $M=\bX(U)$ and $F^{-1}(0)$ have an {\it $A_k$-type contact} at  $p_0=\bX(u_0)$
if the function $g(u) = F \circ \bX(u) $ has the $A_k$ singularity at $u_0.$
By definition, if $M=\bX(U)$ and $F^{-1}(0)$ have an {\it $A_k$-type contact} at  $p_0=\bX(u_0),$
then these have a corank $1$ contact.
For a regular curve $\bgamma :I\lon AdS^{n+1},$ we say that $\bgamma(I)$ and $F^{-1}(0)$ have
a {\it contact of order $k$} if $F\circ \bgamma$ has the $A_k$ singularity at $s_0.$
We have the following simple proposition:
\begin{Pro}
For any $p_0=\bX(u_0)$ and $\blambda _0=\mathbb{LF}_{\kappa _N(\sbn^T)_i(p_0,\sbxi_0)}(u_0),$
there exists an integer $r $ with $1\leq r\leq s$ such that
$M=\bX(U)$ and $TLC_{\sblambda _0}(M)_{p_0}$ have 
corank $r$ contact at $p_0=\bX(u_0).$
\end{Pro}
By Proposition 4.1, $M=\bX (U)$ and the osculating hypersphere $TLC_{\sblambda _0}(M)_{p_0}$
have corank $s$ contact at an $(\bn^T(u_0),\bxi_0)$-umbilical point.
Therefore the AdS-lightlike ridge point is not an $(\bn^T(u_0),\bxi_0)$-umbilical point.
\par
By the general theory of unfoldings of function germs, the discriminant set ${\mathcal D}_F$
is non-singular at the origin if and only if the function $f=F|\R^k\times \{0\}$ has the $A_1$-type singularity (i.e., the
Morse-type singularity).  
Therefore we have the following proposition:

\begin{Pro} With the same notations as in the previous proposition, 
the lightlike hypersurface $\mathbb{LH}_M$ is non-singular at $\blambda _0=\mathbb{LH}_M((p_0,\bxi_0),\mu_0)$ if and only if
$M=\bX(U)$ and $TLC_{\sblambda _0}(M)_{p_0}$ have the $A_1$-type contact at $p_0=\bX(u_0).$
\end{Pro}

\par
We now consider the general contact of $M=\bX(U)$ with AdS-lightcones as an application of  the theory of contact for submanifolds in Montaldi\cite{Montaldi}.
Let $X_i$ and $Y_i,$ $i=1,2,$ be submanifolds of $\R^n$ with
${\rm dim}\, X_1={\rm dim}\, X_2$ and ${\rm dim}\, Y_1={\rm dim}\, Y_2.$
We say that {\it the contact of} $X_1$ and $Y_1$ at $y_1$ is same type as
{\it the contact of} $X_2$ and $Y_2$ at $y_2$ if there is a diffeomorphism germ
$\Phi :(\R^n,y_1)\longrightarrow (\R^n,y_2)$
such that $\Phi (X_1)=X_2$ and $\Phi (Y_1)=Y_2.$
In this case we write
$
K(X_1,Y_1;y_1)=K(X_2,Y_2;y_2).
$
Since this definition of contact is local, we can replace $\R^n$
by arbitrary $n$-manifold.
Montaldi gives in \cite{Montaldi} the following characterization of
contact by using ${\mathcal K}$-equivalence.
\begin{Th} Let $X_i$ and $Y_i,$ $i=1,2,$ be submanifolds of $\R^n$ with
${\rm dim}\, X_1={\rm dim}\, X_2$
and ${\rm dim}\, Y_1={\rm dim}\, Y_2.$
Let $g_i:(X_i,x_i)\longrightarrow
 (\R^n, y_i)$ be immersion germs and $f_i:(\R^n,y_i)\longrightarrow (\R^p,0)$
be submersion germs with $(Y_i,y_i)=(f_i^{-1}(0),y_i).$
Then
$$K(X_1,Y_1;y_1)=K(X_2,Y_2;y_2)$$ if and only if $ f_1\circ g_1$ and
$f_2\circ g_2$ are ${\mathcal K}$-equivalent.
\end{Th}
\par
On the other hand, we now return to the review on the theory of Legendrian singularities.
We introduce a natural equivalence relation among Legendrian submanifold germs.
Let
$F, G : ({\mathbb R}^k\times {\mathbb R}^{n+1},\bo ) \lon  ({\mathbb R},0)$  be
Morse families of hypersurfaces. Then we say
that  $\mathscr{L}_F(\Sigma _*(F))$  and  $\mathscr{L}_G(\Sigma _*(G))$ are
{\it Legendrian equivalent} if there exists a contact diffeomorphism germ
$H : (PT^*{\mathbb R}^{n+1},z) \lon  (PT^*{\mathbb R}^{n+1},z')$  such that  $H$
preserves fibers of  $\pi$ and that $H(\mathscr{L}_F(\Sigma _*(F))) = \mathscr{L}_G(\Sigma _*(G))$, where $z=\mathscr{L}_F(0), z'=\mathscr{L}_G(0).$
By using the Legendrian equivalence, we can define the notion of Legendrian stability for
Legendrian submanifold germs by the ordinary way (see, \cite[Part III]{Arnold1}).
We can interpret the Legendrian equivalence  by using the notion of
generating families.
We denote ${\cal E}_k$ as the local ring of function germs $({\mathbb
R}^k,\bo )\lon
{\mathbb R}$
with the unique maximal ideal $\mathfrak{M}_k=\{h\in {\cal E}_k\ |\ h(0)=0\ \}.$
Let  $F,G : ({\mathbb R}^k\times {\mathbb R}^{n+1},\bo )
\lon  ({\mathbb R},\bo )$  be function germs. We say that  $F$  and  $G$  are
$ P$-${\cal K}$-{\it equivalent} if there exists a diffeomorphism germ
$\Psi  : ({\mathbb R}^k\times {\mathbb R}^{n+1},\bo ) \longrightarrow
({\mathbb R}^k\times
{\mathbb R}^{n+1},\bo )$
of the form  $\Psi (x,u) = (\psi _1(q,x), \psi _{2}(x))$  for
$(q,x) \in  ({\mathbb R}^k\times {\mathbb R}^{n+1},\bo )$  such that
$\Psi ^{*}(\langle F\rangle_{ {\cal E}_{k+n+1}}) = \langle G\rangle_{{ \cal
E}_{k+n+1}}$.
Here  $\Psi ^{*} : {\cal E}_{k+n+1} \lon
{\cal E}_{k+n+1}$  is the pull back ${\mathbb R}$-algebra isomorphism defined by
$\Psi ^{*}(h) =  h\circ \Psi$. We say that $F$ is an {\it infinitesimally ${\cal K}$-versal deformation of}
$f = F\vert
{\mathbb R}^k\times\{ \bo \}$ if
$${\cal E}_k =
T_e({\cal K})(f)
+ \left\langle
\frac{\partial F}{\partial x_1}\vert
{\mathbb R}^k\times\{ \bo \},
\dots ,
\frac{\partial F}{\partial x_{n+1}}\vert
{\mathbb R}^k\times\{ \bo \} \right\rangle_{\mathbb R},
$$
where
$$
T_e({\cal K})(f) =
\left\langle
\frac{\partial f}{\partial q_1}, \dots,
\frac{\partial f}{\partial q_k}, f\right\rangle_{{\cal E}_k}$$
(see, \cite{martine}).
The main result in the theory of Legendrian singularities (\cite{Arnold1}, \S 20.8 and \cite{Zak}, THEOREM 2) is the
following:
\begin{Th}
 Let
$F, G : ({\mathbb R}^k\times {\mathbb R}^{n+1},\bo ) \lon  ({\mathbb R},0)$  be
Morse families of hypersurfaces.
Then we have the following assertions:
\par\noindent
 {\rm (1)} $\mathscr{L} _{F}(\Sigma _*(F))$ and $\mathscr{L} _{G}(\Sigma _*(G))$  are Legendrian equivalent if and
only if
 $F$ and $G$  are  $P$-${\cal K}$-equivalent, 
\par\noindent
{\rm (2)} $\mathscr{L} _{F}(\Sigma _*(F))$   is Legendrian stable if and only if  $F$  is an infinitesimally
${\cal K}$-versal deformation of  $f=F\vert {\mathbb R}^k\times \{\bo \}.$

\end{Th}
\par
Since $F$ and $G$ are function germs on the common space germ $({\mathbb
R}^k\times {\mathbb R}^{n+1},\bo ),$
we do not need the notion of stably $P$-${\cal K}$-equivalences in this
situation \cite[page 27]{Zak}.
For any map germ $f:({\mathbb R}^k,\bo)\lon ({\mathbb R}^p,\bo),$
we define {\it the local ring of} $f$ by
$Q_r(f)={\cal E}_k/(f^*(\mathfrak{M}_p){\cal E}_n+\mathfrak{M}_k^{r+1}).$
We have the following classification result of Legendrian stable germs (cf.
\cite[Proposition A.4]{Izu2}) which is the key for the purpose in this section.
\begin{Pro} Let
$F, G : ({\mathbb R}^k\times {\mathbb R}^{n+1},\bo ) \lon  ({\mathbb R},0)$  be
Morse families of hypersurfaces and $f=F|\R^k\times\{\bo\},g=G|\R^k\times \{\bo\}$.
Suppose that $\mathscr{L} _F(\Sigma _*(F))$ and $\mathscr{L} _G(\Sigma _*(G))$ are Legendrian stable.
The the following conditions are equivalent\/{\rm :}
\par
{\rm (1)} $(W(\mathscr{L} _F),\bo)$ and $(W(\mathscr{L} _G),\bo )$ are diffeomorphic as set germs,
\par
{\rm (2)} $\mathscr{L} _F(\Sigma _*(F))$ and $\mathscr{L} _G(\Sigma _*(G))$ are Legendrian equivalent,
\par
{\rm (3)} $Q_{n+2}(f)$ and $Q_{n+2}(g)$ are isomorphic as ${\mathbb R}$-algebras.
\end{Pro}
\par
We now describe the contacts of spacelike
submanifolds in $AdS^{n+1}$ with AdS-lightcones.
We denote $Q   (\bX ,u_0)$ as the  local ring of the function
germ $\widetilde{h}_{\sblambda _0}:(U,u_0)
\lon \R,$
where $\blambda _0=\mathbb{LC}_M(u_0,\bxi_0,\mu_0).$
We remark that we can explicitly write the local ring as follows:
$$
Q_{n+2}(\bX ,u_0)=
\frac{C^\infty _{u_0}(U)}{\displaystyle{\langle \langle \bX (u),\blambda _0\rangle +1 \rangle _{C^\infty
_{u_0}(U)}}+\mathfrak{M}_{u_0}(U)^{n+2}},
$$
where $C^\infty _{u_0}(U)$ is the local ring of function germs on $U$ at
$u_0.$
\par
Let $\mathbb{LH}_{M_i}(\bn^T_i) :(N_1(M_i)[\bn^T_i]\times\R,(p_i,\bxi_i,\mu_i))\lon (AdS^{n+1} ,\blambda _i),$
$(i=1,2)$ be two
lightlike hypersurface germs of spacelike submanifold germs 
$\bX _i:(U,u^i)\lon (AdS^{n+1},p_i).$ 
Let
$H_i:(U\times AdS^{n+1},(u^i,\blambda  _i))\lon {\mathbb
R}$ be the AdS-height function germ of $\bX _i.$
Then we have the following theorem:
\begin{Th}
Let $\bX _i:(U,u^i)\lon (AdS^{n+1},p_i),$ $i=1,2,$ be
spacelike submanifold germs such that
the corresponding Legendrian submanifold germs $\mathscr{L}_{H_i}(\Sigma _*(H_i))$
are Legendrian stable. We denote that $\bX_i(U)=M_i.$
Then
the following conditions are equivalent$:$
\par
{\rm (1)} $(\mathbb{LH}_{M_1}(N_1(M_1)[\bn^T_1]\times\R),\blambda _1) $ and
$(\mathbb{LH}_{M_2}(N_1(M_2)[\bn^T_2]\times\R),\blambda_2) $
are diffeomorphic,
\par
{\rm (2)} $(\mathscr{L}_{H_1}(\Sigma _*(H_1)), (u^1,\lambda _1)$ and $(\mathscr{L}_{H_2}(\Sigma _*(H_2)), (u^2,\lambda _2)$
are Legendrian equivalent,
\par
{\rm (3)} $H_1$ and $H_2$ are $P$-${\mathcal K}$-equivalent,
\par
{\rm (4)} $h_{1,\lambda _1}$ and $h_{2,\lambda _2}$ are ${\cal K}$-equivalent,
\par
{\rm (5)} $K(M_1,TLC_{\sblambda _1}(M_1)_{p_1} ,p_1)=K(M_2,
TLC_{\sblambda  _2}(M_2)_{p_2},p_2).$
\par
{\rm (6)} $Q_{n+1}  (\bX _1 ,u^1)$ and $Q_{n+1}  (\bX _2
,u^2)$ are isomorphic as ${\mathbb R}$-algebras.
\end{Th}

\demo
By Proposition 5.6, the conditions (1), (2) and (6) are equivalent.
This condition is also equivalent to that
two generating families
$H_1$ and $H_2$ are $P$-${\cal K}$-equivalent by Theorem 5.3. 
If we denote $h_{i,\sblambda  _i}(u)=
H_i(u,\blambda  _i),$ then
we have $h_{i,\sblambda _i}(u)=\mathfrak{h}_{\sblambda _i}\circ\bX
_i(u).$
By Theorem 5.2, $K(\bX _1(U),LC_{\sblambda _1} ,p_1)=K(\bx _2(U),LC{\lambda _2},p_2)$
if and only if $\widetilde{h}_{1,\sblambda _1}$ and $\widetilde{h}_{2,\sblambda
_2}$ are ${\mathcal K}$-equivalent.
This means that (4) and (5) are equivalent.
By definition, (3) implies (4).
The uniqueness of the infinitesimally $\mathcal{K}$-versal deformation of $h_{i,\sblambda _i}$ (cf., \cite{martine})
leads that  (4) implies (3). This completes the proof.
\enD
\par
For a spacelike embedding germ $\bX:(U,u_0)\lon (AdS^{n+1},p_0)$, we consider
a set germ 
$
(\bX ^{-1}(TLC_{\sblambda_0}(M)_{p_0}),u_0),
$
which is called the {\it AdS-tangent lightcone indicatrix germ} of $\bX ,$
where $\blambda_0 =\mathbb{LH}_M(p_0,\bxi _0,\mu_0)$ and
$\mu_0=-1/\kappa_N(\bn^T)_i(p_0,\bxi_0) (i=1,\dots s).$
We have the following corollary of Theorem 5.7.
\begin{Co}
With the assumptions of Theorem {\rm 5.7}, if the lightlike hypersurface germs
$$
(\mathbb{LH}_{M_1}(N_1(M_1)[\bn^T_1]\times\R),\blambda _1)\ and\
(\mathbb{LH}_{M_2}(N_1(M_2)[\bn^T_2]\times\R),\blambda_2) 
$$ are diffeomorphic as set germs, then AdS-tangent lightcone indicatrix germs
$$(\bX_1 ^{-1}(TLC_{\lambda_1}(M_1))_{p_1},u^1)\ and \
(\bX _2^{-1}(TLC_{\lambda_2}(M_2)_{p_2}),u^2)
$$ are diffeomorphic as set germs.
\end{Co}
\demo
We remark that the tangent lightcone indicatrix germ of $\bX_i$ is the zero
level set of
$h_{i,\lambda _i}.$
Since ${\mathcal K}$-equivalence among function germs preserves the
zero-level sets of
function germs, the assertion follows from Theorem 5.7.
\enD
\par
On the other hand, we consider generic properties of lightlike hypersurfaces along spacelike submanifolds.
Let ${\rm Emb}_{\rm sp}\, (U,AdS^{n+1})$ be the space of spacelike embeddings with the 
$C^\infty $-topology for an open set $U\subset \R^s.$
We consider the function
$
{\mathcal H}:AdS^{n+1}\times AdS^{n+1}\longrightarrow \R
$
again.
We claim that $\mathfrak{h}_{\sblambda}$ is a submersion at $\bx\not=\blambda$ for any $\blambda \in AdS^{n+1}.$
For any $\bX\in {\rm Emb}_{\rm sp}\, (U,AdS^{n+1})$, we have $H=\mathcal{H}\circ (\bX\times1_{AdS^{n+1}}).$ 
We have the $r$-jet extension 
$j^r_1H:U\times AdS^{n+1}\lon J^r(U,\R)$ defined by $j^r_1H(u,\blambda)=j^rh_{\sblambda}(u),$
where $J^k(U,\R)$ is the $k$-jet space of functions on $U.$
We consider the trivialization
$J^r(U,\R) \equiv U\times\R\times J^r(s,1).$
For any submanifold $Q\subset J^r(s,1),$ we denote that $\widetilde{Q}=U\times\R\times Q.$
As an application of \cite[Lemma 6]{Wassermann}, the set
\[
T_Q=\{\bX\in {\rm Emb}_{\rm sp}\, (U,AdS^{n+1})\ |\ j^r_1G\ \mbox{is transversal to}\ \widetilde{Q}\ \}
\]
is a residual set of ${\rm Emb}_{\rm sp}\, (U,AdS^{n+1}).$
Moreover, if $Q$ is a closet subset , then $T_Q$ is open.
It is known \cite{GWPL} that there exists a semi-algebraic set $W^r(s,1)\subset J^k(s,1)$ and a stratification
$\mathcal{A}^r(s,1)$ of $J^k(s,1)\setminus W^r(s,1)$ such that $\lim _{k\mapsto \infty}{\rm cod}\, W^r(s,1)=+\infty.$
The stratification $\mathcal{A}^r(s,1)$ is called the {\it canonical stratification}.
We define a stratification $\mathcal{A}^r(U,\R)$ of $J^r(U,\R)\setminus W^r(U,\R)$ by
\[
U\times(\R\setminus \{0\})\times (J^r(s,1)\setminus W^r(s,1)),\ U\times \{0\}\times \mathcal{A}^r(s,1),
\]
where $W^r(U,\R)=U\times\R\times W^r(s,1).$
In \cite{Wan}, it was shown  that if $j^r_1H(U\times AdS^{n+1})\cap W^r(U,\R)=\emptyset$ and
$j^r_1H$ is transversal to $\mathcal{A}^r(U,\R)$, then the map $\pi |H^{-1}(0):H^{-1}(0)\lon AdS^{n+1}$ is MT-stable map-germ at each point,
where $\pi :U\times\R^{n+1}_1\lon AdS^{n+1}$ is the canonical projection.
Here, a map germ is said to be {\it MT-stable} if the jet extension is transversal to the canonical stratification of the jet space of sufficiently higher order (cf., \cite{GWPL,Mather2}).
The main result of the theory of Topological stability of Mather is that MT-stability implies topological stability.
By Proposition 4.1, the lightlike hypersurface $\mathbb{LH}_M(\bn^T)(N_1(M)[\bn^T]\times\R)$ is the discriminant set of
$H$, which is equal to the critical value set of $\pi|H^{-1}(0).$
Since ${\rm cod}\, W^r(U,\R)> s+n+1$ for sufficiently large $k,$  the
set
\[
\mathcal{O}_1=\{\bX\in {\rm Emb}_{\rm sp}\, (U,AdS^{n+1}_1)\ |\ j^r_1H(U\times AdS^{n+1})\cap W^r(U,\R)=\emptyset\ \}
\]
is a residual set.
It follows that the set
\[
\mathcal{O}=\{\bX\in \mathcal{O}_1\ |\ j^r_1H\ \mbox{is transversal to}\ \mathcal{A}^r(U,\R)\ \}
\]
is a residual set. Therefore, we have the following theorem.

\begin{Th}
There exists a residual set $\mathcal{O}\subset {\rm Emb}_{\rm sp}\, (U,AdS^{n+1})$ such that
for any $\bX\in \mathcal{O}$, the germ of the lightlike hypersurface  $\mathbb{LH}_M(\bn^T)(N_1(M)[\bn^T]\times\R)$ at any point
is a germ of the critical value set of an MT-stable map germ.
\end{Th}
In the case when $n\leq 5,$ by the classification results of the $\mathcal{K}$-equivalence among function germs, the canonical stratification $\mathcal{A}^k (s,1)$ is given by the finite collection of
the $\mathcal{K}$-orbits. Moreover, if $j^r_1H$ is transversal to the $\mathcal{K}$-orbit of $j^rh_{\sblambda_0}(u_0)$
for sufficiently large $r,$ then $H$ is an infinitesimally $\mathcal{K}$-versal deformation of $h_{\sblambda}$ at
$(u_0,\blambda _0)$ \cite{martine}. By Theorem 5.5, we have the following theorem. 
\begin{Th} Suppose that $n\leq 5.$ Then there exists a residual set ${\mathcal O}\subset {\rm Emb}_{\rm sp}\, (U,AdS^{n+1})$
such that
for any
$\bX \in {\mathcal O},$ the germ of the lightlike
hypersurface  $\mathbb{LH}_M(\bn^T)(N_1(M)[\bn^T]\times\R)$ at any point is the germ of the wave front set of a stable
Legendrian submanifold germ $\mathscr{L}_H(\Sigma_*(H)).$
\end{Th}

\section{Spacelike submanifolds with codimension two}
In the case when $s=n-1,$ 
$N_1(M)[\bn^T]$ is a double covering of $M=\bX(U).$ 
We can construct  a
spacelike unit normal section $\bn^S(u)\in N_p(M)$ by
\[
\bn ^S(u)=\frac{\bX(u)\wedge \bn^T (u )\wedge\bX _{u_1}(u )\wedge\cdots \wedge
\bX _{u_{n-1}}(u )}{\|\bX(u)\wedge \bn^T (u )\wedge\bX _{u_1}(u )\wedge\cdots
\wedge\bX _{u_{n-1}}(u )\|}.
\]
Then $\bsigma ^\pm (u)=(\bX(u),\pm\bn^S(u))$ are sections of 
$N_1(M)[\bn^T]$.
We call $(\bn ^T,\bn^S)$ a {\it adopted normal frame\/} along $M=\bX(U).$ 
The vectors
$\bn^T (u)\pm \bn^S(u)$ are null. Since
$\{\bX_{u_1}(u),\dots ,\bX_{u_{n-1}}(u)\}$ is a basis of $T_pM,$ the
system $ \{\bX(u), \bn^T (u),\bn ^S(u),\bX_{u_1}(u),\dots
,\bX_{u_{n-1}}(u)\} $ provides a basis for $T_p\R^{n+2}_2$ such that
$\{\bn ^T(u),\bn ^S(u)\}$ is a pseudo-orthonormal frame of the normal timelike plane $N_p(M)\cap T_pAdS^{n+1}$
in $AdS^{n+1}.$
\begin{Lem}\label{parallel}
Given two adopted unit timelike normal sections $\bn^T(u),
\bar{\bn}^T(u)\in N_p(M),$  the corresponding nullcone normal
sections $\bn^T(u)\pm\bn^S(u), \bar{\bn}^T(u)\pm\bar{\bn}^S(u)$ are
parallel.
\end{Lem}
\demo We consider the orientation and the timelike orientation on
the normal space $N_p(M)$ induced by the orientation and the
timelike orientation of $\R^{n+1}_1$ and $\{\bX_{u_1}(u),\dots
,\bX_{u_{n-1}}(u)\}.$ By the construction, both the
pseudo-orthogonal basis $\{\bn ^T(u),\bn^S(u)\}$ and
$\{\bar{\bn}^T(u),\bar{\bn}^S(u)\}$ of $N_p(M)\cap T_pAdS^n$ correspond to the
same orientation and the same timelike orientation on $N_p(M)\cap T_pAdS^n.$
Since both of $\bn^T(u)$ and $\bar{\bn}^T(u)$ are adopted
and $\bn^T(u)\pm\bn^S(u), \bar{\bn}^T(u)\pm\bar{\bn}^S(u)$ are
null in the Lorentz plane $N_p(M)\cap T_pAdS^{n+1}$, $\bn^T(u)\pm\bn^S(u)$ and  $\bar{\bn}^T(u)\pm\bar{\bn}^S(u)$ are
parallel. This completes the proof. \enD
Therefore, the null cone Gauss images of $M=\bX(U)$ with respect to $(\bn^T,\bn^S)$ are given by
$
\mathbb{NG}(\bn^T,\pm\bn^S)(u)=\bn^T(u)\pm\bn^S(u).
$
Since $N_p(M)[\bn^T]$ is a spacelike line in $N_p(M),$ we have $\bxi=\bn^S(u)$ or $\bxi =-\bn^S(u)$ for any $(\bX(u),\bxi)\in N_1(M)[\bn^T].$ 
We denote that $\kappa _i(\bn^T,\pm\bn^S) (u) $, $i=1,\dots , n-1$
instead of $\kappa _N(\bn^T)_i(p,\bxi)$, $i=1,\dots , n-1$ for $\bxi =\pm\bn^S(u)$ 
and $p=\bX(u).$
Then we have the following decomposition of the lightlike hypersurface along $M=\bX(U)$:
\[
\mathbb{LH}_M(N_1(M)[\bn^T]\times\R)=\mathbb{LH}_M^+(U\times\R)\cup \mathbb{LH}_M^-(U\times\R),
\]
where
\[
\mathbb{LH}_M^\pm (u,\mu)= \bX(u)+\mu (\bn^T\pm\bn^S)(u).
\]
In this case the critical value of $\mathbb{LH}_M^\pm $ is the point where
$\kappa_i(\bn^T,\pm\bn^S)(p)\not= 0$ and
\[
\blambda^\pm =\bX(u)+\frac{1}{\kappa_i(\bn^T,\pm\bn^S)(u)}(\bn^T\pm\bn^S)(u).
\]
For each $i=1,\dots ,n-1,$  we have a mapping
$
\mathbb{LE}_{\kappa_i(\sbn^T,\pm\sbn^S)} :O_i\lon AdS^{n+1}
$
defined by
\[
\mathbb{LE}_{\kappa_i(\sbn^T,\pm\sbn^S)}(u)=\bX(u)+\frac{1}{\kappa_i(\bn^T,\pm\bn^S)(u)}(\bn^T\pm\bn^S)(u),
\]
where $O_i=\{u\in U\ |\ \kappa_i(\bn^T,\pm\bn^S)(u)\not= 0\ \}.$
Then we define
\[
\mathbb{LE}^\pm_M=\bigcup\left\{\mathbb{LE}_{\kappa_i(\sbn^T,\pm\sbn^S)}(u)\ |\ u\in U\ \mbox{such\ that}\ \kappa_i(\bn^T,\pm\bn^S)(u)\not= 0,i=1,\dots ,n-1. \right\}.
\]
By the above arguments, we know that $\mathbb{LE}^\pm_{M}$ is nothing but the AdS-lightlike focal set of $M=\bX(U)$
relative to $\bn^T.$
However, we call it the {\it AdS-lightlike evolute} of $M=\bX(U)$ in the case when ${\rm codim}\, M=2.$
For any $p_0=\bX(u_0),$ we have the tangent AdS-lightcones $TLC^\pm_{\sblambda_0}(M)_{p_0},$
where $\blambda ^\pm_0=\mathbb{LE}_{\kappa_i(\sbn^T,\pm\sbn^S)}(u_0).$
In the general codimension case, it depends on the choice of $\bn^S$, so that there are infinitely many tangent AdS-lightcones of $M$ at
$p_0=\bX(u_0)$.
However, the codimension two case, $\bn^S$ is uniquely determined by $\bn^T$.
Therefore, we have two tangent AdS-lightcones of $M$ at $p_0=\bX(u_0).$
In this case, each one of the tangent AdS-lightcones $TLC^\pm _{\sblambda _0}(M)_{p_0}$ is called
 an {\it osculating AdS-lightcone} of $M=\bx (U)$.
The AdS-lightlike evolutes are the loci of the vertices of osculating AdS-lightcones of $M.$
Analogous to the case of hypersurfaces in Euclidean space, we call $\mathbb{LE}^\pm_M$ the {\it AdS-lightlike evolute} of $M=\bX(U).$

\par
We now consider the low dimensions cases in the following subsections.
\subsection{Spacelike curves in $AdS^3$}
We consider spacelike curves in $AdS^3$ as a simplest case of the codimension two case.
Let
$\bgamma:I\longrightarrow AdS^3$ be a unit speed spacelike curve with $\langle \bgamma ''(s),\gamma ''(s)\rangle\not= -1$, where $I$ is an open interval.
We denote that $C=\gamma (I).$
Then we define $\bt(s)= \bgamma '(s)$ and call $\bt(s)$ a {\it unit tangent vector\/ }
of  $\bgamma$ at $s.$
The {\it curvature \/} of $\bgamma$ at $s$ is defined to be
$\kappa _g(s)=\sqrt{|\langle \bgamma (s)-\bgamma ''(s),\bgamma (s)-\bgamma''(s)\rangle | }.$
Since $\kappa _g(s)\not= 0,$ the {\it unit principal normal vector\/} $\bn(s)$
of the curve $\bgamma$ at $s$ is defined by 
$\bgamma''(s)-\bgamma (s)=\kappa _g(s)\bn(s).$
We denote that $\delta{(\bgamma(s))}={\rm sign}(\bn(s))$.
The unit vector $\bb(s)=\bgamma (s)\wedge \bt(s)\wedge \bn(s)$ is called a {\it unit binormal vector}
of the curve $\bgamma $ at $s.$
Since $\bgamma (s)$ is timelike and $\bt(s)$ is spacelike, we have $\langle\bb(s),\bb(s)\rangle =-\delta{(\bgamma(s))}$ and
${\rm sign}\,({\bgamma}'(s))=1$
Then the following Frenet-Serret type formulae hold:
$$ \left\{
\begin{array}{l}
\bgamma '(s)=\bt (s) \\
\bt'{(s)}= \kappa _g(s)\bn(s)+\bgamma (s), \\
\bn'{(s)}=\delta(\bgamma(s))(-\kappa_g(s) 
\bt(s)+ \tau_g{(s)}\bb(s)), \\
\bb'{(s)}= \delta (\bgamma(s))\tau_g{(s)}\bn(s),
\end{array}
\right.
$$
where $\tau_g{(s)}=\langle \bb'(s),\bn(s)\rangle$ is the torsion of the curve $\bgamma$ at $s$.
In \cite{Chen-Han} the singularities of lightlike hypersurfaces along spacelike curves in $AdS^3$ are classified.
We define an invariant $\sigma ^\pm(s)=\kappa _g'(s)\mp\kappa _g(s)\tau _g(s).$
Then we have the following theorem.
\begin{Th}[\cite{Chen-Han}]
Let
$\bgamma:I\longrightarrow AdS^3$ be a unit speed spacelike curve with $\langle \bgamma ''(s),\gamma ''(s)\rangle\not= -1$.
Then we have the following\/{\rm :}
\par\noindent
{\rm (1)} The AdS-lightlike evolute of $\gamma$ is
$$
\mathbb{LE}_C=\left\{\blambda =\bgamma (s)+\frac{1}{\delta (\bgamma (s))\kappa _g(s)}(\bn(s)\pm \bb(s))\right\}.
$$
\par\noindent 
{\rm (2)} The germ of $\mathbb{LH}_M$ at $\blambda _0=\mathbb{LE}_{\kappa _g}(s_0)$ is diffeomorphic to
the cuspidal edge $C(2,3)\times\R$ if and only if $\sigma ^\pm(s_0)\not= 0.$
\par\noindent
{\rm (3)} The germ of $\mathbb{LH}_M$ at $\blambda _0=\mathbb{LE}_{\kappa _g}(s_0)$ is diffeomorphic to
the swallowtail $SW$ if and only if $\sigma ^\pm(s_0)= 0$ and $(\sigma^\pm) '(s_0)\not= 0.$
\end{Th}
\par
Here, $C(2,3)\times \R=\{(x_1,x_2)\ | x_1^2-x_2^3=0\}\times \R$ is the {\it cuspidal edge\/} and
$SW=\{(x_1,x_2,x_3)\ | x_1 =3u^2+u^2v,x_2=4u^3+2uv,x_3=v, (u,v)\in (\R^2,0)\}$
is the {\it swallowtail\/}.
\par
It has been also shown the following geometric characterizations of the singularities of
AdS-lightlike hypersurfaces in \cite{Chen-Han}.
\begin{Pro} Let
$\bgamma:I\longrightarrow AdS^3$ be a unit speed spacelike curve with $\langle \bgamma ''(s),\gamma ''(s)\rangle\not= -1$.
Then we have the followings{\rm :}
\par\noindent
{\rm (1)} $\kappa _g(s_0)\not= 0$ and $\sigma ^\pm (s_0)\not=0$ if and only if $C=\bgamma (I)$ and
the osculating AdS-lightcone $TLC^\pm_{\sblambda _0}(C)_{p_0}$ have contact of order $2.$
\par\noindent
{\rm (2)} $\kappa _g(s_0)\not= 0$, $\sigma ^\pm (s_0)=0$ and $(\sigma ^\pm)'(s_0)\not= 0$ if and only if $C=\bgamma (I)$ and
the osculating AdS-lightcone $TLC^\pm_{\sblambda _0}(C)_{p_0}$ have contact of order $3.$
\end{Pro}
\subsection{Spacelike surfaces in $AdS^4$}
We consider spacelike surfaces in $AdS^4$ here.
Let $\bX:U\lon AdS^4$ be a spacelike embedding from an open subset $U\subset \R^2.$
As a corollary of Theorem 5.10, we have the following generic classification theorem.
We say that two map germs $f,g:(\R^n,0)\lon (\R^p,0)$ are {\it $\mathcal{A}$-equivalent} if
there exists diffeomorphism germs $\phi :(\R^n,0)\lon (\R^n,0)$ and $\psi :(\R^p,0)\lon (\R^p,0)$
such that $f\circ \phi =\psi\circ g.$

\begin{Th}
There exists an open dense subset ${\cal O}\subset {\rm Emb}_{\rm sp}\, (U,AdS^4)$
such that
for any
$\bX \in {\cal O},$ the germ of the corresponding lightlike hypersurfaces
$\mathbb{LH}^\pm _M$ at any point $(u_0,\mu_0)\in U\times \R$
is ${\mathcal A}$-equivalent to one of the map germs $A_k$ $(1\leq k\leq 4)$ or
$D_4^\pm :$ where, $A_k,\ D^\pm_4$-map germ $f:(\R^3,0)\lon (\R^4,0)$ are given by
\par
$A_1;\  f(u_1,u_2,u_3)=(u_1,u_2,u_3,0),$
\par
$A_2; \  f(u_1,u_2,u_3)=(3u_1^2,2u_1^3,u_2,u_3),$
\par
$A_3;\   f(u_1,u_2,u_3)=(4u^3_1+2u_1u_2,3u^4_1+u_2u_1^2,u_2,u_3),$
\par
$A_4;\ 
f(u_1,u_2,u_3)=(5u_1^4+3u_2u_1^2+2u_1u_3,4u_1^5+2u_2u_1^3+u_3u_1^2,u_2,u_3),
$
\par
$D^+_4;\ 
f(u_1,u_2,u_3)=(2(u_1^3+u_2^3)+u_1u_2u_3,3u_1^2+u_2u_3,3u_2^2+u_1u_3,u_3),$
\par
$D^-_4;\ 
\displaystyle{f(u_1,u_2,u_3)=\left(\left(\frac{u_1^3}{3}-u_1u_2^2\right)+(u_1^2+u_2^2)u_3,u_2^2-u_1^2-2u_1u_3,2(u_1u
_2-u_2u_3),u_3\right).}$
\end{Th}
\demo
By Theorems 5.5 and 5.10,  the AdS-height function  $H$ on $M$ is a
${\mathcal K}$-versal  deformation of $h_{\lambda _0}$ at each 
$(u_0,\blambda _0)\in U\times AdS^4.$
Therefore we can apply the classification of ${\mathcal K}$-versal
deformations $F(x,y,\blambda)$ of function germs up to
$4$-parameters \cite{Arnold1}.
For any $F(x,y,\blambda),$ we have
\[
\Sigma _*(F)=\Bigl\{(x,y,\blambda)\in (\R^2\times\R^4,0)\ \Bigm|\ F(x,y,\blambda)=\frac{\partial F}{\partial x}(x,y,\blambda)=
\frac{\partial F}{\partial y}(x,y,\blambda)=0\Bigr\}.
\]
The normal forms are given by
 \begin{eqnarray*}
 A_k;&{}& \!\!\!\!\!\!\! F(x,y,\blambda )=x^{k+1}\pm y^2+\lambda _1+\lambda _2x+\cdots +\lambda
_{k-1}x^{k-1},\ 1\leq k\leq 4, \\
D^+_4;&{}& \!\!\!\!\!\!\! F(x,y,\blambda )=x^3+y^3+\lambda _1+\lambda _2x+\lambda _3y+\lambda _4xy, \\
D^-_4;&{}& \!\!\!\!\!\!\! F(x,y,\blambda )= \frac{x^3}{3}-xy^2+\lambda _1+\lambda _2x+\lambda _3y+\lambda
_4(x^2+y^2).
\end{eqnarray*}
For example, if we consider the germ given by
$$
F(x,y,\blambda )=x^3+y^3+\lambda _1+\lambda _2x+\lambda _3y+\lambda _4xy.
$$
Then we get
$$
\Sigma _*(F)=\{(x,y,2(x^3+y^3)+\lambda _4xy,-3x^2-\lambda _4y,-3y^2-\lambda_4x,\lambda _4)\ |\ (x,y,\lambda _4)\in \R^3
\}.
$$
Therefore the corresponding Legendrian map germ is
$$
f(u_1,u_2,u_3)=(2(u_1^2+u_2^2)+u_1u_2u_3,3u_1^2+u_2u_3,3u_2^2+u_1u_3,u_3)\quad (D^+_4) .
$$
The other cases follow from similar arguments, so that we omit the details.
\enD
As a corollary of the above theorem, we have the following generic local classification of
AdS-lightlike evolutes along spacelike surfaces.
We define that $C(2,3,4)=\{(u_1^2,u_1^3,u_1^4)\ |\ u_1\in \R\}$ which is called
a {\it $(2,3,4)$-cusp}.
We also define that
$C(BF)=\{(10u_1^3+3u_2u_1,5u_1^4+u_2u_1^2,6u_1^5+u_2u_1^3,u_2)\ |\ (u_1,u_2)\in \R^2\}$.
We call $C(BF)$ a {\it C-butterfly} (i.e., the critical value set of the butterfly).
Finally we define that where $C(2,3,4,5)=\{(u_1^2,u_1^3,u_1^4,u_1^5)\ |\ u_1\in \R\}$ which is called
a {\it $(2,3,4,5)$-cusp}.

\begin{Co} There exists an open dense subset ${\cal O}\subset {\rm Emb}_{\rm sp}\, (U,AdS^4)$
such that
for any
$\bX \in {\cal O},$ the germ of the corresponding AdS-lightlike evolute
$\mathbb{LE}^\pm _M$ at any point $(u_0,\mu_0)\in U\times \R$ is diffeomorphic to one of 
the following set germs at the origin in $\R^4${\rm:}
\par\noindent
$A_2; \ \{(0,0)\}\times\R^2,$
\par\noindent
$A_3;\  C(2,3,4)\times \R,$
\par\noindent
$A_4;\ C(BF),
%\{(5u_1^4+u_2u_1^2,6u_1^5+u_2u_1^3,u_2,10u_1^3+3u_2u_1)|\ (u_1,u_2)\in \R^2\},
$
\par\noindent
$D^+_4;\ 
\{(2(u_1^3+u_2^3)+u_1u_2u_3,3u_1^2+u_2u_3,3u_2^2+u_1u_3,u_3)\ |\ u_3^2=36u_1u_2\},$
\par\noindent
$D^-_4;\ 
\displaystyle{\left\{\left(\left(\frac{u_1^3}{3}-u_1u_2^2\right)+(u_1^2+u_2^2)u_3,u_2^2-u_1^2-2u_1u_3,2(u_1u
_2-u_2u_3),u_3\right)\ \Bigl|\ u_3^2=u_1^2+u_2^2\right\}.}$
\end{Co} 
\demo
For $A_3$, we can calculate the Jacobi matrix of the normal form $f$ in Theorem 6.4:
\[
J_f=\left(
\begin{array}{ccc}
12u^2_1+2u_1 & 2u_1 & 0 \\
12u^3_1+2u_1u_2 & u^2_1 & 0 \\
0 & 1 &0 \\
0 & 0 & 1
\end{array}
\right),
\]
so that ${\rm rank}\, J_f <3$ if and only if $6u^2_1+u_2=0.$
Thus, the critical value set of $f$ is
$$C(f)=\{(-8u^3_1,-3u^4_1,-6u^2_1,u_3)\ |\ (u_1,u_3)\in \R^2\}.$$
It is $C(2,3,4)\times \R.$
By the similar calculation, we can show that the germ of $A_4$ is diffeomorphic to $C(BF).$
For $D^+_4,$ we can calculate the Jacobi matrix o the normal form $f$:
\[
J_f=\left(
\begin{array}{ccc}
6u^2_1+u_2u_3 & 6u_2^2+u_1u_3,u_1u_2 & 0 \\
6u_1 & u_3 & u_2 \\
u_3 & 6u_2 &u_1 \\
0 & 0 & 1
\end{array}
\right).
\]
Therefore, ${\rm rank}\,J_f<3$ if and only if
\[
\left|\begin{array}{cc}
6u^2_1+u_2u_3 & 6u_2^2+u_1u_3,u_1u_2 \\
6u_1 & u_3  \\
\end{array}
\right|=
\left|
\begin{array}{cc}
6u^2_1+u_2u_3 & 6u_2^2+u_1u_3,u_1u_2 \\
u_3 & 6u_2\\
\end{array}
\right|=
\left|
\begin{array}{cc}
6u_1 & u_3  \\
u_3 & 6u_2 \\
\end{array}
\right|=0,
\]
which is equivalent to the condition that 
$u_3^2=36u_1u_2.$
For $D^-_4,$ by the similar calculation to the above, we have the condition that
$u_3^2=u_1^2+u_2^2.$
This completes the proof.
\enD
\par
In the list of the above corollary, 
one of the projections of the image of $D^+_4$ into $\R^3$ is $PU=\{(3u_1^2+u_2u_3,3u_2^2+u_1u_3,u_3)\ |\ u_3^2=36u_1u_2\}$ 
which is called a {\it purse} and one of the projection of the image of $D^-_4$ is
called a {\it pyramid} given by 
\[
PY=\left\{\left(u_2^2-u_1^2-2u_1u_3,2(u_1u
_2-u_2u_3),u_3\right)\ \Bigl|\ u_3^2=u_1^2+u_2^2\right\}.
\]
We can draw these pictures as follows:
We denote that $\widetilde{PU}$ as the set of $D^+_4$ and $\widetilde{PY}$ as the set of $D^-_4$
in the list of Corollary 6.5, respectively.
We also have a classification of the singularities of $\mathbb{LE}^\pm _M$ as a corollary of Theorem 6.4
and Corollary 6.5.
The set of singularities
 of $\mathbb{LE}^\pm _M$ is denoted by $\Sigma (\mathbb{LE}^\pm _M).$ 
\begin{Co}
There exists an open dense subset ${\cal O}\subset {\rm Emb}_{\rm sp}\, (U,AdS^4)$
such that
for any
$\bX \in {\cal O},$ the germ of the pair 
$(\mathbb{LE}^\pm _M,\Sigma(\mathbb{LE}^\pm _M))$ at any point $\blambda_0
\in \Sigma(\mathbb{LE}^\pm_M)$ is diffeomorphic to one of 
the following pairs of set germs at the origin in $\R^4${\rm:}
\par\noindent
$A_3;\ (C(2,3,4)\times \R, \{(0,0,0)\}\times \R),$
\par\noindent
$A_4;\ (C(BF),C(2,3,4,5)),
$
\par\noindent
$D^+_4;\ (\widetilde{PU},\Sigma (\widetilde{PU})),
$
where
$\Sigma (\widetilde{PU})=\{(5u^3_3/108,u^2_3/4,u^2_3/4,u_3)\ |\ u_3\in \R\},$
\par\noindent
$D^-_4;\ (\widetilde{PY},\Sigma (\widetilde{PY})),
$
where
\par\noindent
$
\Sigma (\widetilde{PY}))=\{\left(4u^3_3/3,-3u^2_3,0,u_3\right)\ |\ u_3\in \R\}
\cup \\
\hfill \{(4u^3_3/3,3u^2_3/2,-3\sqrt{3}u^2_3/2,u_3)\ |\ u_3\in \R\}\cup 
\{(4u^3_3/3,3u^2_3/2,3\sqrt{3}u^2_3/2,u_3)\ |\ u_3\in \R\}.
$
\end{Co}
\demo
For $A_3$, in order to detect the singularities, we calculate the Jacobi matrix
of $f(u_1,u_2)=(u_1^2,u_1^3,u_1^4,u_2).$ Then we have
\[
J_f=\left(
\begin{array}{cc}
2u_1 & 0 \\
3u_1^2 & 0 \\
4u_1^3 & 0 \\
0 & 1
\end{array}
\right),
\]
so that ${\rm rank}\, J_f<2$ if and only if $u_1=0.$
This means that $\Sigma (C(2,3,4)\times\R)=\{(0,0,0)\}\times \R.$
Similar to the case $A_3,$ we calculate the Jacobi matrix
of $f(u_1,u_2)=(10u_1^3+3u_2u_1,5u_1^4+u_2u_1^2,6u_1^5+u_2u_1^3,u_2)$.
Then we have
\[
J_f=\left(
\begin{array}{cc}
30u_1^2+3u_2 & 3u_1 \\
20u_1^3+2u_2u_2 & u_1^2 \\
30u_1^4+3u_2u_1^2 & u_1^3 \\
0 & 1
\end{array}
\right),
\]
so that ${\rm rank}\, J_f<2$ if and only if $u_2=-10u_1^2.$
Therefore, we have
\[
f(u_1,-10u_1^2)=(-20u_1^3,-5u_1^4,-4u_1^5,-10u_1^2).
\]
This means that $\Sigma (C(BF))=C(2,3,4,5).$
For $D^+_4,$ we consider the following parameter transformation:
\[
u_1+u_2=\frac{u_3}{3}\cosh\phi,\ u_1-u_2=\frac{u_3}{3}\sinh\phi.
\]
Then $\widetilde{PU}$ is parametrized by
\[
f(\phi,u_3)=\left(\frac{4}{6^3}u_3^3\left(\cosh 3\phi+\frac{1}{4}\right),\frac{u_3^2}{6^2}\left(3e^{2\phi}+6e^{-\phi}\right),
\frac{u_3^2}{6^2}\left(3e^{-2\phi}+e^\phi\right),u_3\right).\]
Thus , the Jacobi matrix is 
\[
J_f=\left(
\begin{array}{cc}
\displaystyle{\frac{1}{18}u_3^3\sinh3\phi} & \displaystyle{\frac{1}{18}u_3^2(\cosh3\phi +\frac{1}{4})}\\
\displaystyle{\frac{u_3^2}{6}(e^{2\phi}-e^{-\phi})} & \displaystyle{\frac{1}{6}(e^{2\phi}+2e^{-\phi})} \\
\displaystyle{\frac{u_3^2}{6}(-e^{-2\phi}+e^\phi)} & \displaystyle{\frac{1}{6}(e^{-2\phi}+2e^\phi)}\\
0 & 1
\end{array}
\right),
\]
so that ${\rm rank}\, J_f<2$ if and only if $\phi=0$ or $u_3=0.$
This means that $u_1=u_3/6,\ u_2=u_3/6.$
Therefore, we have
\[\Sigma (\widetilde{PU})=\left\{\left(\frac{5}{108}u^3_3,\frac{1}{4}u^2_3,\frac{1}{4}u^2_3,u_3\right)\ \Bigl|\ u_3\in \R\right\}.
\]
For $D^-_4,$ we also consider the following parameter transformation:
\[
u_1=u_3\cos\theta,\ u_2=u_3\sin\theta.
\]
Then $\widetilde{PY}$ is parametrized by
\[
f(\theta,u_3)=\left(\frac{1}{3}u_3^3\left(\cos3\theta+1\right),-u_3^2(\cos2\theta+2\cos\theta),
2u_3^2(\sin2\theta -2\sin\theta),u_3\right).\]
Thus , the Jacobi matrix is 
\[
J_f=\left(
\begin{array}{cc}
-u_3^3\sin3\theta & u_3^2(\cos3\theta +1) \\
2u_3^2(\sin2\theta +\sin\theta) & -2u_3(\cos2\theta+2\cos\theta) \\
2u_3^2(\cos2\theta -\cos\theta) & 2u_3(\sin2\theta-2\sin\theta) \\
0 & 1
\end{array}
\right),
\]
so that ${\rm rank}\, J_f<2$ if and only if $\theta=0,\ 2\pi/3, 5\pi/3$ or $u_3=0.$
If $\theta=0,$ then we have $u_1=u_3,\ u_2=0,$ so that
we have $\{\left(4u^3_3/3,-3u^2_3,0,u_3\right)| u_3\in \R\}.$
If $\theta =2\pi/3,$ then we have
$u_1=-u_3/2,\ u_2=\sqrt{3}u_3/2,$ so that we have
$\{(4u^3_3/3,3u^2_3/2,-3\sqrt{3}u^2_3/2,u_3| u_3\in \R\}.$
Finally, if $\theta =5\pi/3,$ then we have
$u=-u_3/2,\ u_2=-\sqrt{3}u_3/2,$ so that we have
$\{(4u^3_3/3,3u^2_3/2,3\sqrt{3}u^2_3/2,u_3)| u_3\in \R\}.$
It follows that $\Sigma (\widetilde{PY}))$ is the union of these three set germs.
This completes the proof.
\enD
\par
We can interpret the geometric meaning of the above normal forms of the singularities of
lightlike hypersurfaces as the following theorem shows.

\begin{Th}
Let $\mathcal{O}\subset {\rm Emb}_{\rm sp}\, (U,AdS^4)$ be the open dense subset given in
the above theorem. For any $\bX \in {\cal O},$ the germ of the corresponding lightlike hypersurfaces
$(\mathbb{LH}^\pm _M,\blambda _0)$ at a point $(u_0,\mu_0)\in U\times \R$
is characterized as follows{\rm :}
\par\noindent
{\rm (1)} $(A_1)$ if and only if $\blambda _0\in \mathbb{LH}^+_M(U\times \R)\cup \mathbb{LH}^-_M(U\times \R)$ is a non-singular point.
\par\noindent
{\rm (2)} $(A_2)$ if and only if $p_0=\bX(u_0)$ is not a AdS-lightlike ridge point and $\blambda_0
\in \mathbb{LE}^+_M\cup\mathbb{LE}^-_M.$
\par\noindent
{\rm (3)} $(A_3)$ if and only if $p_0=\bX(u_0)$ is an AdS-lightlike $1$-ridge point and $\blambda_0
\in \mathbb{LE}^+_M\cup \mathbb{LE}^-_M.$
\par\noindent
{\rm (4)} $(A_4)$ if and only if $p_0=\bX(u_0)$ is an AdS-lightlike $2$-ridge point, $\blambda_0
\in \mathbb{LE}^+_M\cup \mathbb{LE}^-_M$ and $u_0\in U$ satisfies the following condition{\rm :} For sufficiently small $\varepsilon >0$, there exist two different AdS-lightlike $1$-ridge points $u^1,u^2\in U$
such that $|u_0-u^i|<\varepsilon$, $(i=1,2)$.
\par\noindent
{\rm (5)} $(D_4^+)$ if and only if $\blambda_0
\in \mathbb{LE}^+_M\cup\mathbb{LE}^-_M$, $p_0=\bX(u_0)$ has corank two contact with the osculating lightcone and
$u_0\in U$ satisfies
the following condition{\rm :} 
\par
{\rm (a)} $\kappa _1(\bn^T,\pm\bn^S)(u^0)= \kappa _2(\bn^T,\pm\bn^S)(u^0)$.
\par
{\rm (b)}
 For sufficiently small $\varepsilon >0$, there exist two different AdS-lightlike $1$-ridge points $u^1,u^2\in U$
such that $|u_0-u^i|<\varepsilon$, $(i=1,2)$ and $\kappa _1(\bn^T,\pm\bn^S)(u^i)\not= \kappa _2(\bn^T,\pm\bn^S)(u^i)$.
\par\noindent
{\rm (6)} $(D_4^-)$ if and only if $\blambda_0
\in \mathbb{LE}_M$, $p_0=\bX(u_0)$ has corank two contact with the the osculating lightcone and
$u_0\in U$ satisfies
the following condition{\rm :} 
\par
{\rm (a)} $\kappa _1(\bn^T,\pm\bn^S)(u^0)= \kappa _2(\bn^T,\pm\bn^S)(u^0)$.
\par
{\rm (b)} For sufficiently small $\varepsilon >0$, there exist three different AdS-lightlike $1$-ridge points $u^i\in U$, $(i=1,2,3)$
such that $|u_0-u^i|<\varepsilon$ and
$\kappa _1(\bn^T,\pm\bn^S)(u^i)\not=\kappa _2(\bn^T,\pm \bn^S)(u^i)$.
\end{Th}

\demo
By the normal form $(A_1)$ and $(A_2)$ in Theorem 6.4, the assertions (1) and (2) are trivial.
For the normal form $(A_3)$, the AdS-height function has the $A_3$ singularity at $p_0=\bX(u_0)$, so that
it is an AdS-lightlike $1$-ridge point.
\par
Here, 
we give a remark on the classification of $\mathcal{K}$-simple singularities of function germs.
In the list of the classification, we say that a class of singularities $L$ is {\it adjacent} to a class $K$
(notation: $K\leftarrow L$) if every function germ $f\in L$ can be deformed into a function of $K$ by an
arbitrarily small perturbation.
For the class of $A_k,D^\pm_k$ of the $\mathcal{K}$-classification are adjacent to each other as follows \cite[Page 243]{Arnold1}:
\[
\begin{array}{cccccccc}
A_1 &\leftarrow & A_2 &\leftarrow & A_3 &\leftarrow & A_4 &\leftarrow \\
{} & {} & {} & {} & \uparrow & {} & \uparrow & {} \\
{} & {} & {} & {} & D^\pm_4 &\leftarrow & D^\pm_5 &\leftarrow  \\
\end{array}
\]
By the normal form $(A_4)$, the singularities of the AdS-lightcone evolute is a $(2,3,4,5)$-cusp.
Therefore, two singular loci approach to the $(2,3,4,5)$-cusp point. 
Since $A_4$ is adjacent to $A_3$, such the singular loci 
consist of AdS-lightlike $1$-ridge points except the origin.
Thus, the assertion (4) holds.
For $(D^+_4)$, the singularities of the AdS-lightcone evolute is $\Sigma (\widetilde{PU})$.
By the normal form of the generating family, the corank of the AdS-height function at $u^0\in U$ is
two, so that it is an $(\bn^T(u_0),\pm \bn^S(u_0))$-umbilical point. Therefore, we have
$\kappa _1(\bn^T,\pm\bn^S)(u^0)= \kappa _2(\bn^T,\pm\bn^S)(u^0)$.
By the normal form of $(D^+_4)$ in Corollary 6.6, two singular loci approach to the origin
from both the positive and the negative side of the parameter $u_3$.
Since $D^+_4$ is adjacent to $A_3$, such the singular loci consist of AdS-lightlike $1$-ridge points
except at the origin.
Of course, an AdS-lightlike $1$-ridge point is not an $(\bn^T,\pm \bn^S)$-umbilical point.
So that two nullcone principal curvatures are different at an AdS-lightlike $1$-ridge point.
For $(D^-_4)$, we have the assertion by the similar arguments to the case $(D^+_4).$
This completes the proof.
\enD

\section{Spacelike curves in $AdS^4$}
In this section we consider lightlike hypersurfaces along spacelike curves in $AdS^4$ as
the simplest case of higher codimension.
%\section{Spacelike curves in anti de Sitter $4$-space}
%In \S 6 we investigated the spacelike submanifolds with codimension two, and we have a classification of the singularities of the lightlike hypersurfaces in $\R^4_1.$ In this section we consider the higher codimensional case in $\R^4_1$, that is spacelike curves in Minkowski $4$-space as a special case as the previous results.
%\par
Let $\bgamma:I\lon AdS^4$ be a spacelike curve with $\| \bgamma''(s)\|\neq-1.$ In this case we write $C=\bgamma (I)$ instead of $M=\bgamma (I).$ Since $\|\bgamma'(s)\|>0,$ we can reparameterize it by the arc-length s. So we have the unit tangent vector $\bt(s)=\bgamma'(s)$ of $\bgamma(s).$ Moreover we have two unit normal vectors $\displaystyle{\bn_1(s)=\frac{\bgamma''(s)-\bgamma(s)}{\|\bgamma''(s)-\bgamma(s)\|}},$ $\displaystyle{\bn_2(s)=\frac{\bn'_1(s)+\delta \kappa _1(s)\bt(s)}{\| \bn'_1(s)+\delta \kappa _1(s)\bt(s)\|}}$ under the conditions that $\kappa_1(s)=\|\bgamma''(s)-\bgamma(s)\|\neq0,$ $\kappa _2(s)=\| \bn'_1(s)+\delta k_1(s)\bt(s)\|\neq0,$ where $\delta_i={\rm sign}(\bn_i(s))$ and ${\rm sign}(\bn_i(s))$ is the signature of $\bn_i(s)$ $(i=1,2,3).$ Then we have another unit normal vector field $\bn_3(s)$ defined by $\bn_3(s)=\bgamma(s)\wedge\bt(s)\wedge\bn_1(s)\wedge\bn_2(s).$ Therefore we can construct a pseudo-orthogonal frame $\{\bgamma(s), \bt(s), \bn_1(s), \bn_2(s), \bn_3(s)\},$ which satisfies the {\it Frenet-Serret type formulae}:
 $$\left\{
       \begin{array}{ll}
          \bgamma'(s)=\bt(s), \\
          \bt'(s)=\bgamma(s)+\kappa_1(s)\bn_1(s),  \\
          \bn_1'(s)=-\delta_1\kappa_1(s)\bt(s)+\kappa_2(s)\bn_2(s),\\
          \bn_2'(s)=\delta_3\kappa_2(s)\bn_1(s)+\kappa_3(s)\bn_3(s),\\
          \bn_3'(s)=\delta_1\kappa_3(s)n_2(s),  
       \end{array}
     \right.$$   
where $\kappa_2(s)= \delta_2\langle \bn_1'(s),\bn_2(s)\rangle$ and $\kappa_3(s)=\delta_3\langle \bn_2'(s),\bn_3(s)\rangle.$ 
%\par\noindent
Since $\bgamma(s)$ is timelike and $\bt(s)$ is spacelike, we distinguish the following three cases:
\par
\smallskip
Case 1: $\bn_1(s)$ is timelike, that is, $\delta_1=-1$ and $\delta_2=\delta_3=1.$
\par
Case 2: $\bn_2(s)$ is timelike, that is, $\delta_2=-1$ and $\delta_1=\delta_3=1.$
\par
Case 3: $\bn_3(s)$ is timelike, that is, $\delta_3=-1$ and $\delta_1=\delta_2=1.$
\par\noindent
We consider the lightlike hypersurface along $C,$ and calculate the anti-de Sitter height  function on $C$ which is useful for the study the singularities of lightlike hypersurfaces in the each case.
%\par\noindent
\subsection{Case 1}
Suppose that $\bn_1(s)$ is timelike. In this case we adopt  $\bn^T(s)=\bn_1(s)$ and denote that $\bb_1(s)=\bn_2(s), \bb_2(s)=\bn_3(s).$ Then we have the pseudo-orthogonal frame 
\[
\{ \bgamma(s), \bt(s), \bn^T(s), \bb_1(s), \bb_2(s)\},
\]
$\delta_1=-1$ and $\delta_2=\delta_3=1$, which satisfies the following Frenet-Serret type formulae:
$$\left\{
       \begin{array}{ll}
          \bgamma'(s)=\bt(s), \\
          \bt'(s)=\bgamma(s)+\kappa_1(s)\bn^T(s),  \\
          {\bn^T}'(s)=\kappa_1(s)\bt(s)+\kappa_2(s)\bb_1(s),\\
          \bb_1'(s)=\kappa_2(s)\bn^T(s)+\kappa_3(s)\bb_2(s),\\
          \bb_2'(s)=-\kappa_3(s)\bb_1(s). 
       \end{array}
     \right.$$
Since $N_1(C)[\bn^T]$ is parametrized by 
\[
N_1(C)[\bn^T]=\{(\bgamma (s),\bxi)\in \bgamma ^*T\R^5_1\ |\ \bxi =\cos\theta\bb_1(s)+\sin\theta\bb_2(s)\in N_{\sbgamma(s)}(C),\ s\in I \},
\]  
the nullcone Gauss image of $N_1(C)_{p}[\bn^T]$ is given by
\[
\mathbb{NG}(\bn^T)(s,\theta)=\bn^T(s)+\cos\theta\bb_1(s)+\sin\theta\bb_2(s).
\]
Then we have the lightlike hypersurface along $C$ 
$$
\mathbb{LH}_C((s,\theta),\mu)=\bgamma(s)+\mu(\bn^T(s)+\cos\theta\bb_1(s)+\sin\theta\bb_2(s))=
\bgamma (s)+\mu\mathbb{NG}(\bn^T)(s,\theta).
$$
 We remark that the image of this lightlike hypersurface along $C$ is independent of the choice of the future directed timelike normal vector field $\bn^T$ by Corollary 4.3. 
\par\noindent
Now we investigate the anti- de Sitter height functions $H: I\times AdS^{4}\lon \R$ on a spacelike curve  $C=\bgamma(I)$ defined by
$$
 H(p,\blambda)=H(s,\blambda)=\langle \bgamma(s),\blambda\rangle+1,
$$
where $p=\bgamma(s).$ For any fixed $\blambda _0\in AdS^{4},$ we write $h(p)=H_{\lambda _0}(p)=H(p,\blambda _0).$
\par\noindent
By Proposition 4.1, the discriminant set of the
anti-de Sitter height function $H$ is given by
\[
{\mathcal D}_{H}=\mathbb{LH}_C(N_1(C)[\bn^T]\times\R)=\Bigl\{\blambda =\bgamma(s)+\mu\mathbb{NG}(s,\theta)\Bigm|\ \theta\in [0,2\pi), s\in I,\ \mu\in \R\ \Bigr\},
\]
which is the image of the lightlike hypersurface along $C.$ 
We also calculate that $h''(p)=\langle \bgamma''(s), \lambda_0\rangle=-\mu\kappa_1-1.$ Then $h''(p)=0$ if and only if $\mu=-1/\kappa_1(s).$ 
 It means that a singular point of the lightlike hypersurface is a point $\blambda _0=\bgamma(s_0)+\mu_0\mathbb{NG}(\theta _0,s_0)$
for $\mu_0 =-1/\kappa_1(s_0).$ 
Therefore, the lightlike focal surface is
\[
\mathbb{LF}_C=\Bigl\{\blambda =\bgamma(s)-\frac{1}{\kappa _1(s)}\mathbb{NG}(\bn^T)
(s,\theta)\Bigm| \ s\in I,\ \theta\in [0,2\pi)\ \Bigr\}.
\]

\par
Moreover, if we calculate the third, 4th and 5th derivatives of $h(s),$ we have the following proposition.
\par
\begin{Pro}
Let $C$ be a spacelike curve and
$H: C\times(AdS^{4}\setminus C)\to\R$
the anti-de Sitter height function on $C.$
Suppose that $p_0=\bgamma(s_0)\not=\blambda _0.$ Then we have the followings{\rm :}
\par\noindent
{\rm (1)}
$h(p_0)=h'(p_0)=0$ if and only if
there exist $\theta _0\in [0,2\pi)$ and $\mu\in
\R\setminus \{0\}$ such that 
$$
\bgamma (s_0)-\blambda _0 =\mu\mathbb{NG}(\bn^T)(s_0,\theta_0).
$$ 
\par\noindent
{\rm (2)}
$h(p_0)=h'(p_0)=h''(p_0)=0$ if and only if
there exists $\theta _0\in [0,2\pi)$ such that 
$$
\bgamma (s_0)-\blambda _0=-\frac{1}{\kappa_1(s_0)}\mathbb{NG}(\bn^T)(s_0,\theta _0).
$$
\par\noindent
{\rm(3)}
$h(p_0)=h'(p_0)=h''(p_0)=h'''(p_0)=0$ if and only if there exists $\theta _0\in [0,2\pi)$ such that 
$$
\bgamma (s_0)-\blambda _0=-\frac{1}{\kappa_1(s_0)}\mathbb{NG}(\bn^T)(s_0,\theta _0)
$$
and $\kappa_1'(s_0)-\cos\theta _0 \kappa_1(s_0)\kappa_2(s_0)=0,$
so that we can write $\theta _0=\theta (s_0).$
\par\noindent
{\rm(4)}
$h(p_0)=h'(p_0)=h''(p_0)=h'''(p_0)=h^{(4)}(p_0)=0$ if and only if there exists $\theta _0=\theta (s_0)\in [0,2\pi)$ such that 
$$
\bgamma (s_0)-\blambda _0=-\frac{1}{\kappa_1(s_0)}\mathbb{NG}(\bn^T)(s_0,\theta (s_0)),
$$
$\kappa_1'(s_0)-\cos\theta (s_0) \kappa_1(s_0)\kappa_2(s_0)=0$ and 
$(2\kappa_1'(s_0)\kappa_2(s_0)+\kappa_1(s_0)\kappa_2'(s_0))\cos\theta (s_0)-\kappa_1''(s_0)-\kappa_1(s_0)\kappa_2^2(s_0)+\kappa_1(s_0)\kappa_2(s_0)\kappa_3(s_0)\sin\theta (s_0)=0.$
\par\noindent
{\rm(5)}
$h(p_0)=h'(p_0)=h''(p_0)=h'''(p_0)=h^{(4)}(p_0)=h^{(5)}(p_0)=0$ if and only if there exists $\theta _0=\theta (s_0)\in [0,2\pi)$ such that 
$$
\bgamma (s_0)-\blambda _0=-\frac{1}{\kappa_1(s_0)}\mathbb{NG}(\bn^T)(s_0,\theta (s_0)),
$$
$\kappa_1'(s_0)-\cos\theta (s_0) \kappa_1(s_0)\kappa_2(s_0)=0$,
$(2\kappa_1'(s_0)\kappa_2(s_0)+\kappa_1(s_0)\kappa_2'(s_0))\cos\theta (s_0)-\kappa_1''(s_0)-\kappa_1(s_0)\kappa_2^2(s_0)+\kappa_1(s_0)\kappa_2(s_0)\kappa_3(s_0)\sin\theta (s_0)=0$
and
$((2\kappa_1'(s_0)\kappa_2(s_0)+\kappa_1(s_0)\kappa_2'(s_0))\cos\theta(s_0)-\kappa_1''(s_0)-\kappa_1(s_0)\kappa_2^2(s_0)+\kappa_1(s_0)\kappa_2(s_0)\kappa_3(s_0)\sin\theta(s_0))'=0.$
\end{Pro}
\par
Taking account of the above proposition, we denote that
$\rho _1(s,\theta)=\kappa_1'(s)-\cos\theta  \kappa_1(s)\kappa_2(s)$
and
$\eta _1(s,\theta)=(2\kappa_1'(s)\kappa_2(s)+\kappa_1(s)\kappa_2'(s))\cos\theta -\kappa_1''(s)-\kappa_1(s)\kappa_2^2(s)+\kappa_1(s)\kappa_2(s)\kappa_3(s)\sin\theta ,$
which might be important invariants of $C=\bgamma (I).$
Then we can show that
$\rho_1(s,\theta)=\eta_1(s,\theta)=0$ if and only if $\rho _1(s,\theta)=\sigma _1(s)=0$,
where
\[
\sigma _1(s)=\left[\kappa _1\kappa _2(\kappa _1''+\kappa _1\kappa _2^2)
-\kappa _1'(2\kappa _1'\kappa _2+\kappa _1\kappa '_2)\mp
\kappa _1\kappa _2\kappa _3\sqrt{(\kappa _1\kappa _2)^2-(\kappa '_1)^2}\right](s).
\]

\par
\subsection{Case 2}
Suppose that $\bn_2(s)$ is timelike. Then we adopt $\bn^T(s)=\bn_2(s)$ and denote that $\bb_1(s)=\bn_1(s), \bb_2(s)=\bn_3(s).$ We have  a pseudo-orthogonal frame $\{ \bgamma(s), \bt(s), \bn^T(s), \bb_1(s), \bb_2(s)\}$, $\delta_2=-1$ and $\delta_1=\delta_3=1,$ which satisfies the following Frenet-Serret type formulae:
$$\left\{
       \begin{array}{ll}
          \bgamma'(s)=\bt(s), \\
          \bt'(s)=\bgamma(s)+\kappa_1(s)\bb_1(s),  \\
          \bb_1'(s)=-\kappa_1(s)\bt(s)+\kappa_2(s)\bn^T(s),\\
          {\bn^T}'(s)=\kappa_2(s)\bb_1(s)+\kappa_3(s)\bb_2(s),\\
          \bb_2'(s)=\kappa_3(s)\bn^T(s),  
       \end{array}
     \right.$$
Here, $N_1(C)[\bn^T]$ is parametrized by 
\[
N_1(C)[\bn^T]=\{(\bgamma (s),\bxi)\in \bgamma ^*T\R^5_1\ |\ \bxi =\cos\theta\bb_1(s)+\sin\theta\bb_2(s)\in N_{\sbgamma(s)}(C),\ s\in I \},
\]  so that we have the lightlike hypersurface along $C=\bgamma (I)$:
$$
\mathbb{LH}_C((s,\theta),t)=\bgamma(s)+\mu\mathbb{NG}(\bn^T)(s,\theta).
$$
\par
We consider the anti-de Sitter height function $H: I\times AdS^{4}\lon \R$ on a spacelike curve  $C=\bgamma(I)$. Under the similar notations to the case 1), we have the following proposition:
\begin{Pro}
Let $C$ be a spacelike curve and
$H: C\times(AdS^{4}\setminus C)\to\R$
the anti-de Sitter height function on $C.$
Suppose that $p_0\not=\blambda _0.$ Then we have the following$:$
\par\noindent
{\rm (1)}
$h(p_0)=h'(p_0)=0$ if and only if
there exist $\theta _0 \in [0,2\pi)$ and $\mu\in
\R\setminus \{0\}$ such that 
$$
\bgamma (s_0)-\blambda _0 =\mu\mathbb{NG}(\bn^T)(s_0,\theta_0).
$$ 
\par\noindent
{\rm (2)}
$h(p_0)=h'(p_0)=h''(p_0)=0$ if and only if
there exists $\theta _0 \in [0,2\pi)$ such that 
$$
\bgamma(s_0)-\blambda _0=\frac{1}{\kappa_1(s_0)\cos\theta_0}\mathbb{NG}(\bn^T)(s_0,\theta _0).
$$
\par\noindent
{\rm(3)}
$h(p_0)=h'(p_0)=h''(p_0)=h'''(p_0)=0$ if and only if there exists $\theta _0 \in [0,2\pi)$ such that 
$$
\bgamma(s_0)-\blambda _0=\frac{1}{\kappa_1(s_0)\cos\theta_0}\mathbb{NG}(\bn^T)(s_0,\theta _0)
$$
and $\kappa_1'(s_0)\cos\theta_0-\kappa_1(s_0)\kappa_2(s_0)=0,$
so that we can write $\theta _0=\theta (s_0).$
\par\noindent
{\rm(4)}
$h(p_0)=h'(p_0)=h''(p_0)=h'''(p_0)=h^{(4)}(p_0)=0$ if and only if there exists $\theta _0 =\theta (s_0)\in [0,2\pi)$ such that 
$$
\bgamma(s_0)-\blambda _0=\frac{1}{\kappa_1(s_0)\cos\theta (s_0)}\mathbb{NG}(\bn^T)(s_0,\theta (s_0)),
$$
$\kappa_1'(s_0)\cos\theta (s_0)-\kappa_1(s_0)\kappa_2(s_0)=0$ and
$(\kappa_1''(s_0)+\kappa_1(s_0)\kappa_2^2(s_0))\cos\theta (s_0)-2\kappa_1'(s_0)\kappa_2(s_0)-\kappa_1(s_0)\kappa_2'(s_0)+\kappa_1(s_0)\kappa_2(s_0)\kappa_3(s_0)\sin\theta (s_0)=0$.
\par\noindent
{\rm(5)}
$h(p_0)=h'(p_0)=h''(p_0)=h'''(p_0)=h^{(4)}(p_0)=h^{(5)}(p_0)=0$ if and only if there exists $\theta _0 =\theta (s_0)\in [0,2\pi)$ such that 
$$
\bgamma(s_0)-\blambda _0=\frac{1}{\kappa_1(s_0)\cos\theta (s_0)}\mathbb{NG}(\bn^T)(s_0,\theta (s_0)),
$$
$\kappa_1'(s_0)\cos\theta (s_0)-\kappa_1(s_0)\kappa_2(s_0)=0$,
$(\kappa_1''(s_0)+\kappa_1(s_0)\kappa_2^2(s_0))\cos\theta (s_0)-2\kappa_1'(s_0)\kappa_2(s_0)-\kappa_1(s_0)\kappa_2'(s_0)+\kappa_1(s_0)\kappa_2(s_0)\kappa_3(s_0)\sin\theta (s_0)=0$ 
and $((\kappa_1''(s_0)+\kappa_1(s_0)\kappa_2^2(s_0))\cos\theta (s_0)-2\kappa_1'(s_0)\kappa_2(s_0)-\kappa_1(s_0)\kappa_2'(s_0)+\kappa_1(s_0)\kappa_2(s_0)\kappa_3(s_0)\sin\theta (s_0))'=0.$
\end{Pro}
\par\noindent
The above proposition asserts that the discriminant set of the Lorentzian distance-squared function $G$ is given by
\[
{\mathcal D}_{H}=\mathbb{LH}_C(N_1(C)[\bn^T]\times\R)=\Bigl\{\blambda =\bgamma(s)-\mu\mathbb{NG}(\bn^T)
(s,\theta)\Bigm| \ s\in I,\ \theta\in [0,2\pi), \mu\in \R\ \Bigr\}.
\]
Moreover, the lightlike focal surface is
\[
\mathbb{LF}_C=\Bigl\{\blambda =\bgamma(s)-\frac{1}{\kappa _1(s)\cos\theta}\mathbb{NG}(\bn^T)
(s,\theta)\Bigm| \ s\in I,\ \theta\in [0,2\pi)\ \Bigr\}.
\]
\par
Here, we also denote that
$
\rho_2(s,\theta)=\kappa '_1(s)\cos\theta -\kappa _1(s)\kappa _2(s)$ and
\[
\eta _2(s,\theta)=(\kappa ''_1(s)+\kappa _1(s)\kappa ^2_2(s))\cos\theta -2\kappa '_1(s)\kappa _2(s)-\kappa _1(s)\kappa '_2(s)+\kappa _1(s)\kappa _2(s)\kappa _3(s)\sin\theta.
\]
We can also show that
$\rho_2(s,\theta)=\eta_2(s,\theta)=0$ if and only if $\rho _2(s,\theta)=\sigma _2(s)=0$,
where
\[
\sigma _2(s)=\left[\kappa _1\kappa _2(\kappa _1''+\kappa _1\kappa _2^2)
-\kappa _1'(2\kappa _1'\kappa _2+\kappa _1\kappa '_2)\pm
\kappa _1\kappa _2\kappa _3\sqrt{-(\kappa _1\kappa _2)^2+(\kappa '_1)^2}\right](s).
\]

%which might be important invariants of $C=\bgamma (I).$
\par
\subsection{Case 3}
Suppose that $\bn_3(s)$ is timelike. Then we adopt $\bn^T(s)=\bn_3(s)$ and denote that 
$\bb_1(s)=\bn_1(s), \bb_2(s)=\bn_2(s).$ We have  a pseudo-orthogonal frame $\{ \bgamma(s), \bt(s), \bn^T(s), \bb_1(s), \bb_2(s)\}$ and $\delta_3=-1$ and $\delta_1=\delta_2=1,$which satisfies the following Frenet-Serret type formulae:
$$\left\{
       \begin{array}{ll}
          \bgamma'(s)=\bt(s), \\
          \bt'(s)=\bgamma(s)+\kappa_1(s)\bb_1(s),  \\
          \bb_1'(s)=-\kappa_1(s)\bt(s)+\kappa_2(s)\bb_2(s),\\
          \bb_2'(s)=-\kappa_2(s)\bb_1(s)+\kappa_3(s)\bn^T(s),\\
          {\bn^T}'(s)=\kappa_3(s)\bb_2(s),  
       \end{array}
     \right.$$
Here, $N_1(C)[\bn^T]$ is parametrized by 
\[
N_1(C)[\bn^T]=\{(\bgamma (s),\bxi)\in \bgamma ^*T\R^5_1\ |\ \bxi =\cos\theta\bb_1(s)+\sin\theta\bb_2(s)\in N_{\sbgamma(s)}(C),\ s\in I \},
\]  
so that we have the lightlike hypersurface along $C$:
$$
\mathbb{LH}_C((s,\theta),t)=\bgamma(s)+t\mathbb{NG}(\bn^T)(s,\theta).
$$
We investigate the anti-de Sitter function on a spacelike curve  $C=\bgamma(I)$ 
By the calculations similar to the cases 1 and 2, we have the following proposition:
\par
\begin{Pro}
Let $C$ be a spacelike curve and
$H: C\times(AdS^4\setminus C)\to\R$
the anti-de Sitter function on $C=\bgamma (I).$
Suppose that $p_0\not=\blambda _0.$ Then we have the following$:$
\par\noindent
{\rm (1)}
$h(p_0)=h'(p_0)=0$ if and only if
there exist $\theta _0\in [0,2\pi)$ and $\mu\in
\R\setminus \{0\}$ such that 
$$
\bgamma (s_0)-\blambda _0 =\mu\mathbb{NG}(\bn^T)(s_0,\theta _0).
$$ 
\par\noindent
{\rm (2)}
$h(p_0)=h'(p_0)=h''(p_0)=0$ if and only if
there exists $\theta _0\in [0,s\pi)$ such that 
$$
\bgamma (s_0)-\blambda _0=\frac{1}{\kappa_1(s_0)\cos\theta_0}\mathbb{NG}(\bn^T)(s_0,\theta _0).
$$
\par\noindent
{\rm(3)}
$h(p_0)=h'(p_0)=h''(p_0)=h'''(p_0)=0$ if and only if there exists $\theta _0\in [0,s\pi)$ such that 
$$
\bgamma (s_0)-\blambda _0=\frac{1}{\kappa_1(s_0)\cos\theta_0}\mathbb{NG}(\bn^T)(s_0,\theta _0)
$$
and $\kappa_1'(s_0)\cos\theta _0+\kappa_1(s_0)\kappa_2(s_0)\sin\theta _0=0,$
so that we can write $\theta _0=\theta (s_0).$
\par\noindent
{\rm(4)}
$h(p_0)=h'(p_0)=h''(p_0)=h'''(p_0)=h^{(4)}(p_0)=0$ if and only if there exists $\theta _0=\theta (s_0)\in [0,2\pi)$ such that 
$$
\bgamma (s_0)-\blambda _0=\frac{1}{\kappa_1(s_0)\cos\theta (s_0)}\mathbb{NG}(\bn^T)(s_0,\theta (s_0)),
$$
$\kappa_1'(s_0)\cos\theta (s_0)+\kappa_1(s_0)\kappa_2(s_0)\sin\theta (s_0)=0.$
 and  $(2\kappa_1'(s_0)\kappa_2(s_0)+\kappa_1(s_0)\kappa_2'(s_0))\sin\theta (s_0)+(\kappa_1''(s_0)-\kappa_1(s_0)\kappa_2^2(s_0))\cos\theta (s_0)-\kappa_1(s_0)\kappa_2(s_0)\kappa_3(s_0)=0.$
\par\noindent
{\rm(5)}
$h(p_0)=h'(p_0)=h''(p_0)=h'''(p_0)=h^{(4)}(p_0)=h^{(5)}(p_0)=0$ if and only if there exists $\theta _0=\theta (s_0)\in [0,2\pi)$ such that 
$$
\bgamma (s_0)-\blambda _0=\frac{1}{\kappa_1(s_0)\cos\theta(s_0)}\mathbb{NG}(\bn^T)(s_0,\theta (s_0)),
$$
$\kappa_1'(s_0)\cos\theta (s_0)+\kappa_1(s_0)\kappa_2(s_0)\sin\theta (s_0)=0,$
and $((2\kappa_1'(s_0)\kappa_2(s_0)+\kappa_1(s_0)\kappa_2'(s_0))\sin\theta (s_0))+(\kappa_1''(s_0)-\kappa_1(s_0)\kappa_2^2(s_0))\cos\theta (s_0))-\kappa_1(s_0)\kappa_2(s_0)\kappa_3(s_0))'=0.$
\end{Pro}
The above proposition asserts that the discriminant set of the anti-de Sitter function $H$ is given by
\[
{\mathcal D}_{H}=\mathbb{LH}_C(N_1(C)[\bn^T]\times\R)=\Bigl\{\blambda =\bgamma(s)+t\mathbb{NG}(\bn^T)(s,\theta) \Big|\  s\in I,\ \theta\in [0,2\pi), t\in \R\ \Bigr\}.
\]
Moreover, the lightlike focal surface is
\[
\mathbb{LF}_C=\Bigl\{\blambda =\bgamma(s)-\frac{1}{\kappa _1(s)\cos\theta}\mathbb{NG}(\bn^T)
(s,\theta)\Bigm| \ s\in I,\ \theta\in [0,2\pi)\ \Bigr\}.
\]
\par
Here, we also denote that
$
\rho_3(s,\theta)=\kappa '_1(s)\cos\theta +\kappa _1(s)\kappa _2(s)\sin\theta$ and
\[
\eta _3(s,\theta)=(2\kappa '_1(s)\kappa _2(s)+\kappa _1(s)\kappa '_2(s))\sin\theta +(\kappa ''_1(s)-\kappa _1(s)\kappa ^2_2(s))\cos\theta-\kappa _1(s)\kappa _2(s)\kappa _3(s),
\]
We can also show that
$\rho_3(s,\theta)=\eta _3(s,\theta)=0$ if and only if $\rho _3(s,\theta)=\sigma _3(s)=0$,
where
\[
\sigma _3(s)=\left[\kappa _1\kappa _2(\kappa _1''-\kappa _1\kappa _2^2)
-\kappa _1'(2\kappa _1'\kappa _2+\kappa _1\kappa '_2)\mp
\kappa _1\kappa _2\kappa _3\sqrt{(\kappa _1\kappa _2)^2+(\kappa '_1)^2}\right](s).
\]
%which might be important invariants of $C=\bgamma (I).$
\par
We can unify the invariants $\sigma _i(s)$, $(i=1,2,3)$ as follows:
\[
\sigma (s)=\left[\kappa _1\kappa _2(\kappa _1''-\kappa _1\kappa _2^2)
-\kappa _1'(2\kappa _1'\kappa _2+\kappa _1\kappa '_2)\mp\delta _2
\kappa _1\kappa _2\kappa _3\sqrt{\delta _1(\kappa _1\kappa _2)^2+\delta _2(\kappa '_1)^2}\right](s).
\]
\par\noindent
\subsection{Classifications of singularities}
By using the results of the three cases, we classify the singularities of the lightlike hypersurface along $\bgamma$ as an application of the unfolding theory of functions. 
For a function $f(s)$, we say that $f$
has {\it $A_k$-singularity } at $s_0$ if $f^{(p)}(s_0)=0$  for all
$1\leq p\leq k$ and $f^{(k+1)}(s_0)\neq 0.$
 Let $F$ be an  $r$-parameter unfolding of $f$
and $f$ has  $ A_k $-singularity $(k\geq 1)$ at $s_0$. We denote
the $(k-1)$-jet of the partial derivative $\partial F/
\partial x_i$ at $s_0$
 as
$$j^{(k-1)}\left(\frac{\partial F}{\partial x_i}(s,\bx_0)\right)(s_0)
 =\sum\limits _{j=1}^{k-1} \alpha_{ji}(s-s_0)^j,~~  (i=1, \cdots ,r).$$
  If the rank of $k\times r$ matrix  $(\alpha_{0i},\alpha_{ji}) $ is  $k~(k\leq
r)$, then $F$ is called a {\it versal unfolding} of $f$,
    where $ \alpha_{0 i}=\partial F/\partial x_i(s_0,\bx_0)$.
\par
Inspired by the propositions in the previous subsections, we define the following set:
\[
D^\ell _F=\left\{\bx\in  \R^r\mid \exists s\in \R,\ F(s,\bx)=\frac{\partial F}{\partial s}(s,\bx)=\cdots =\frac{\partial ^\ell F}{\partial s^\ell }(s,\bx)= 0 \right\},
\]
which is called a {\it discriminant set of order $\ell $.}
Of course, $D^1_F=D_F$ and $D^2_F$ is the set of singular points of $D_F.$
Therefore, we have the following proposition.
\begin{Pro}
For all the cases, we have
\[
D_G=D^1_G=\mathbb{LH}_C(N_1(C)[\bn^T]\times\R),\ D^2_G=\mathbb{LF}_C\ \mbox{and}\ D^3_G\ 
\mbox{is\ the\ critical\ value\ set\ of}\ \mathbb{LF}_C.
\]
\end{Pro}
In order to understand the geometric properties of the discriminant set of order $\ell$, we introduce
an equivalence relation among the unfoldings of functions.
Let $F$ and $G$ be $r$-parameter unfoldings of $f(s)$ and  $g(s)$, respectively.
We say that $F$ and $G$ are {\it P-$\mathcal{R}$-equivalent} if
there exists a diffeomorphism germ $\Phi :(\R\times\R^r,(s_0,\bx_0))\longrightarrow (\R\times\R^r, (s'_0,\bx'_0))$
of the form $\Phi (s,\bx)=(\Phi _1(s,\bx),\phi (\bx))$ such that $G\circ\Phi =F.$
By straightforward calculations, we have the following proposition.
\begin{Pro} Let $F$ and $G$ be $r$-parameter unfoldings of $f(s)$ and  $g(s)$, respectively. If $F$ and $G$ are P-$\mathcal{R}$-equivalent by a diffeomorphism germ
$\Phi :(\R\times\R^r,(s_0,\bx_0))\longrightarrow (\R\times\R^r, (s'_0,\bx'_0))$ of the form
 $\Phi (s,\bx)=(\Phi _1(s,\bx),\phi (\bx))$, then $\phi (D^\ell _F)=D^\ell _G$ as set germs.
\end{Pro}
    
We have the following classification theorem of versal unfoldings \cite[Page 149, 6.6]{Bru-Gib}.
 \begin{Th}
 Let $F:(\R\times\R^r,(s_0,\bx_0))\longrightarrow \R $ be an
$r$-parameter unfolding of $f$ which has  $A_k$-singularity at
$s_0$. Suppose  $F$ is a versal unfolding of $f$, then $F$ is P-$\mathcal{R}$-equivalent to
one of the following unfoldings:
\par
{\rm (a)} $k=1$ {\rm ;} $\pm s^2+x_1$,
\par
{\rm (b)} $k=2$ {\rm ;} $s^3+x_1+sx_2$,
\par
{\rm (c)} $k=3$ {\rm ;} $\pm s^4+x_1+sx_2+s^2x_3,$
\par
{\rm (d)} $k=4$ {\rm ;} $s^5+x_1+sx_2+s^2x_3+s^3x_4.$
\end{Th}
\par\noindent
We have the following classification result as a corollary of the above theorem.
\begin{Co}
 Let $F:(\R\times
\R^r,(s_0,\bx_0))\longrightarrow \R $ be an
$r$-parameter unfolding of $f$ which has  $A_k$-singularity at
$s_0$. Suppose  $F$ is a versal unfolding of $f$, then we have the following assertions:
\par\noindent
{\rm (a)} If $k=1$, then $ D_F $ is diffeomorphic to $\{0\}\times
\R^{r-1}$ and $D^2_F=\emptyset.$
\par\noindent
{\rm (b)} If  $k=2$, then $D _F $ is diffeomorphic to $C(2,3)\times \R^{r-2},$
$D^2_F$ is diffeomorphic to $\{\bo\}\times\R^{r-2}$
and $D^3_F=\emptyset.$
\par\noindent
{\rm (c)} If $k=3$, then $D _F $ is diffeomorphic to $SW\times \R^{r-3},$ 
 $D^2_F$ is diffeomorphic to $C(2,3,4)\times \R^{r-3},$
$D^3_F$ is diffeomorphic to $\{\bo\}\times \R^{r-3}$ and $D^4_F=\emptyset.$
\par\noindent
{\rm (d)} If $k=4$, then $D _F $ is locally diffeomorphic to $BF\times \R^{r-4},$
$D^2_F$ is diffeomorphic to $C(BF)\times \R^{r-4},$
$D^3_F$ is diffeomorphic to $C(2,3,4,5)\times \R^{r-4}$,
$D^4_F$ is diffeomorphic to $\{\bo\}\times\R^{r-4}$ and $D^5_F=\emptyset.$
\par\noindent
We remark that all of diffeomorphisms in the above assertions are diffeomorphism germs.
\end{Co}
Here, we call $BF=\{(x_1,x_2,x_3.x_4)\mid x_1 = 5u^4+3vu^2+2wu,x_2=4u^5+2vu^3+wu^2,x_3=u,x_4=v\}$ a {\it butterfly}.
We have the following key proposition on $H.$
\begin{Pro}
 If $h(s)$ has $A_k$-singularity $(k=1,2,3,4)$ at $p_0$, then $H$ is a versal unfolding of $h.$
\end{Pro}
\demo
We denote that $\bgamma(s)=(X_{-1}(s),X_0(s),X_1(s),X_2(s),X_3(s))\ {\rm and}\
\blambda =(\lambda_{-1},\lambda _0,\lambda _1,\lambda _2,\lambda _3).$
\par\noindent
By definition, we have
$$
H(s,\blambda )=-\lambda_{-1}X_{-1}(s)-\lambda_0X_0(s)+\lambda_1X_1(s)+\lambda_2X_2(s)+\lambda_3X_3(s)+1.
$$
Thus we have 
$$
\frac{\partial H}{\partial \lambda_i}(s,\blambda)=-X_i (s), {\rm and} \ \frac{\partial^2 H}{\partial s\partial \lambda_i}(s,\blambda)=-X_i' (s), \ {\rm for} (i=-1,0)
$$
$$
\frac{\partial H}{\partial \lambda_i}(s,\blambda)=X_i (s), {\rm and} \ \frac{\partial^2 H}{\partial s\partial \lambda_i}(s,\blambda)=X_i' (s), \ {\rm for} (i=1,2,3)
$$
For a fixed $\blambda_0=(\lambda_{0-1}\lambda_{00}, \lambda_{01}, \lambda_{02}, \lambda_{03}),$ the 3-jet of $\partial H/\partial
  \lambda_i(s,\blambda_0)(i=-1,0,1,2,3)$ at $s_0$ is
  $$j^{(3)}\frac{\partial H}{\partial
  \lambda_i}(s,\blambda_0)(s_0)=-X_i'(s_0)(s-s_0)-\frac{1}{2}X_i''(s_0)(s-s_0)^2-\frac{1}{3}X_i'''(s_0)(s-s_0)^3, \ (i=-1,0).$$
$$j^{(3)}\frac{\partial H}{\partial
  \lambda_i}(s,\blambda_0)(s_0)=X_i'(s_0)(s-s_0)+\frac{1}{2}X_i''(s_0)(s-s_0)^2+\frac{1}{3}X_i'''(s_0)(s-s_0)^3, \ (i=1,2,3).$$
\par\noindent
It is enough to show that the rank of the following matrix A is four,
$$A=\left(
\begin{array}{ccccc}
 -X_{-1}(s_0)& -X_0 (s_0)& X_1 (s_0)& X_2 (s_0)& X_3 (s_0)\\
-X'_{-1}(s_0)& -X_0'(s_0)& X_1'(s_0)& X_2'(s_0)& X_3'(s_0)\\
-X''_{-1}(s_0)& -X_0''(s_0)& X_1''(s_0)& X_2''(s_0)& X_3''(s_0)\\
-X'''_{-1}(s_0)& -X_0'''(s_0)& X_1'''(s_0)& X_2'''(s_0)& X_3'''(s_0)
\end{array}
\right).$$
Here we consider the following matrix B,
$$B=\left(
\begin{array}{ccccc}
 -X_{-1}(s_0)& -X_0 (s_0)& X_1 (s_0)& X_2 (s_0)& X_3 (s_0)\\
-X'_{-1}(s_0)& -X_0'(s_0)& X_1'(s_0)& X_2'(s_0)& X_3'(s_0)\\
-X''_{-1}(s_0)& -X_0''(s_0)& X_1''(s_0)& X_2''(s_0)& X_3''(s_0)\\
-X'''_{-1}(s_0)& -X_0'''(s_0)& X_1'''(s_0)& X_2'''(s_0)& X_3'''(s_0)\\
-n_{-1}(s_0)& -n_{0}(s_0)& n_{1}(s_0)& n_{2}(s_0)& n_{3}(s_0)\\
\end{array}
\right),$$
where we denote $\bn_1(s_0)=(n_{-1}(s_0), n_{0}(s_0), n_{1}(s_0), n_{2}(s_0), n_{3}(s_0)).$
In fact, 
\begin{eqnarray*}
{\rm det}B&=&{\rm det}^{t}(\bgamma(s), \bgamma'(s), \bgamma''(s), \bgamma'''(s),\bn_{1})\\
&=&{\rm det}^{t}(\bgamma(s), \bgamma'(s),\bgamma''(s)-\bgamma(s),\bb, \bn_{1})\\
&=&{\rm det}^{t}(\bgamma(s), \bgamma'(s),\bgamma''(s)-\bgamma(s)-\kappa_{1}\bn_{1},\bb-\kappa_{1}\bn_{1},\bn_{1}),\\
\end{eqnarray*}
where $\bb=\bgamma'''(s)+\frac{\kappa'_{1}(s)}{\kappa_{1}(s)}\bgamma(s)+(\delta_{1}\kappa^2_{1}(s)-1)\bgamma'(s)+(1-\frac{\kappa_{1}'(s)}{\kappa_{1}})\bgamma''(s),$
%\end{eqnarray*}
and $\bgamma(s), \bgamma'(s),\bgamma''(s)-\bgamma(s)-\kappa_{1}\bn_{1},$ $\bb-\kappa_{1}\bn_{1}$ and $\bn_{1}$  are linearly independent each other in all Case 1,2,3, respectively. Therefore we have detB$\neq0$. This means that ${\rm rank A}=4.$
This completes the proof.
\enD

\par\noindent
Finally, we can apply Corollary 8.5 to our condition. Then we have the following theorem:
\begin{Th}
Let $\bgamma:I\lon AdS^4$ be a spacelike curves with $\kappa _1(s)\not= 0.$ %and $\kappa _2(s)\not= 0.$
\par\noindent
{\rm (A)} For the lightlike hypersurfaces $\mathbb{LH}_C((s,\theta),t)$ of $C=\bgamma (I)$ in the Case 1,  
we have the following assertions:
\par\noindent
{\rm (1)} The lightlike hypersurface  $\mathbb{LH}_C(N_1(C)[\bn^T]\times\R )$ is locally diffeomorphic to $C(2,3)\times \R^2$ at
$\blambda_0$ if and only if there exist $\theta _0\in [0,2\pi )$ such that 
$$
p_0-\blambda _0=-\frac{1}{\kappa_1(s_0)}\mathbb{NG}(\bn^T)(s_0,\theta _0),
$$
and $\rho_1(s_0,\theta_0)\neq0.$ In this case, the lightlike focal set $\mathbb{LF}_C$ is non-singular.
\par\noindent
{\rm(2)} The lightlike hypersurface  $\mathbb{LH}_C(N_1(C)[\bn^T]\times\R )$ is locally diffeomorphic to $SW\times \R$ at
$\blambda_0$ if and only if there exist $\theta _0\in [0,2\pi )$ such that 
$$
p_0-\blambda _0=-\frac{1}{\kappa_1(s_0)}\mathbb{NG}(\bn^T)(s_0,\theta _0),
$$
$\rho _1(s_0,\theta_0)=0$ and $\sigma _1(s_0)\neq 0.$
In this case, the lightlike focal set $\mathbb{LF}_C$ is locally diffeomorphic to $C(2,3,4)\times\R$ and the critical value set of $\mathbb{LF}_C$ is a regular curve.

\par\noindent
{\rm(3)} The lightlike hypersurface $\mathbb{LH}_C(N_1(C)[\bn^T]\times\R )$  is locally diffeomorphic to $BF$ at
$\blambda_0$ if and only if there exist $\theta _0\in [0,2\pi )$ such that 
$$
p_0-\blambda _0=-\frac{1}{\kappa_1(s_0)}\mathbb{NG}(\bn^T)(s_0,\theta _0),
$$
$\rho _1(s_0,\theta _0)=0$, $\sigma _1(s_0)=0$ and $\sigma '_1(s_0)\not= 0.$
In this case, the lightlike focal set $\mathbb{LF}_C$ is is locally diffeomorphic to $C(BF)\times\R$ and the critical value set is locally diffeomorphic to the $C(2,3,4,5)$-cusp.

\smallskip
\par\noindent
{\rm (B)} For the lightlike hypersurfaces $\mathbb{LH}_C(N_1(C)[\bn^T]\times\R )$ of $C=\bgamma(I)$ in the Case 2,  
we have the following assertions:
\par\noindent
{\rm (1)} The lightlike hypersurface  $\mathbb{LH}_C(N_1(C)[\bn^T]\times\R )$ is locally diffeomorphic to $C(2,3)\times \R^2$ at
$\blambda_0$ if and only if there exist $\theta _0\in [0,2\pi)$ such that 
$$
p_0-\blambda _0=\frac{1}{\kappa_1(s_0)\cos\theta_0}\mathbb{NG}(\bn^T)(s_0,\theta_0)
$$
and  $\rho_2(s_0,\theta _0)\neq0.$ In this case, the lightlike focal set
$\mathbb{LF}_C$ is non-singular.
\par\noindent
{\rm(2)} The lightlike hypersurface  $\mathbb{LH}_C(N_1(C)[\bn^T]\times\R )$ is locally diffeomorphic to $SW\times \R$ at
$\blambda_0$ if and only if there exist $\theta _0\in [0,2\pi)$ such that 
$$
p_0-\blambda _0=\frac{1}{\kappa_1(s_0)\cos\theta_0}\mathbb{NG}(\bn^T)(s_0,\theta_0),
$$
$\rho _2(s_0,\theta _0)=0$ and 
$\sigma (s_0)\neq0.$ In this case, the lightlike focal set $\mathbb{LF}_C$ is locally
diffeomorphic to $C(2,3,4)\times\R$ and the critical value set of $\mathbb{LF}_C$ is a regular curve.
\par\noindent
{\rm(3)} he lightlike hypersurface  $\mathbb{LH}_C(N_1(C)[\bn^T]\times\R )$ is locally diffeomorphic to $BF$ at
$\blambda_0$ if and only if there exist $\theta _0\in [0,2\pi)$ such that 
$$
p_0-\blambda _0=\frac{1}{\kappa_1(s_0)\cos\theta_0}\mathbb{NG}(\bn^T)(s_0,\theta_0),
$$
$\rho _2(s_0,\theta _0)=0$, 
$\sigma _2(s_0)=0$ and $\sigma '_2(s_0)\neq 0.$ In this case, the lightlike focal set
$\mathbb{LF}_C$ is locally
diffeomorphic to $C(BF)\times\R$ and the critical value set of $\mathbb{LF}_C$ is locally diffeomorphic to the $C(2,3,4,5)$-cusp.

\smallskip
\par\noindent
{\rm (C)} For the lightlike hypersurfaces $\mathbb{LH}_C(N_1(C)[\bn^T]\times\R )$ of $C=\bgamma(I)$ in the Case 3,  
we have the following assertions:
\par\noindent
{\rm (1)} The lightlike hypersurface  $\mathbb{LH}_C(N_1(C)[\bn^T]\times\R )$ is locally diffeomorphic to $C(2,3)\times \R^2$ at
$\blambda_0$ if and only if there exist $\theta _0\in [0,2\pi)$ such that 
$$
p_0-\blambda _0=\frac{1}{\kappa_1(s_0)\cos\theta_0}\mathbb{NG}(\bn^T)(s_0,\theta_0),
$$
and $\rho _3(s_0,\theta _0)\neq0.$ In this case, the lightlike focal set
$\mathbb{LF}_C$ is non-singular.

\par\noindent
{\rm(2)} The lightlike hypersurface  $\mathbb{LH}_C(N_1(C)[\bn^T]\times\R )$ is locally diffeomorphic to  $SW\times \R$ at
$\blambda_0$ if and only if there exist $\theta _0\in [0,2\pi)$ such that 
$$
p_0-\blambda _0=\frac{1}{\kappa_1(s_0)\cos\theta_0}\mathbb{NG}(\bn^T)(s_0,\theta_0),
$$
$\rho _3(s_0,\theta _0)=0$ and $\sigma _3(s_0)\neq 0.$ In this case, the lightlike focal set
$\mathbb{LF}_C$ is locally
diffeomorphic to $C(2,3,4)\times\R$ and the critical value set of $\mathbb{LF}_C$ is a regular curve.
\par\noindent
{\rm(3)} The lightlike hypersurface  $\mathbb{LH}_C(N_1(C)[\bn^T]\times\R )$ is locally diffeomorphic to $BF$ at
$\blambda_0$ if and only if there exist $\theta _0\in [0,2\pi)$ such that 
$$
p_0-\blambda _0=\frac{1}{\kappa_1(s_0)\cos\theta_0}\mathbb{NG}(\bn^T)(s_0,\theta_0),
$$
$\rho _3(s_0,\theta _0)=0$, $\sigma _3(s_0)= 0$ and $\sigma '_3(s_0)\not= 0.$
In this case, the lightlike focal set
$\mathbb{LF}_C$ is locally
diffeomorphic to $C(BF)\times\R$ and the critical value set of $\mathbb{LF}_C$ is locally diffeomorphic to the $C(2,3,4,5)$-cusp.
\end{Th}

{\small
\par\noindent
Shyuichi Izumiya, Department of Mathematics, Hokkaido University, Sapporo 060-0810,Japan
\par\noindent
e-mail:{\tt izumiya@math.sci.hokudai.ac.jp}
}
\end{document}